\documentclass[leqno]{amsart}
\usepackage{amsfonts}
\usepackage{amscd,amsmath,amsopn,amssymb,amsthm,multicol}

\newtheorem{theorem}{Theorem}

\newtheorem{corollary}[theorem]{Corollary}

\newtheorem{definition}[theorem]{Definition}
\newtheorem{example}[theorem]{Example}

\newtheorem{lemma}[theorem]{Lemma}
\newtheorem{notation}[theorem]{Notation}

\newtheorem{proposition}[theorem]{Proposition}
\newtheorem{remark}[theorem]{Remark}

\input{tcilatex}

\begin{document}
\title{Covering $\mathbb{R}$-trees}
\author{V.N.~Berestovski\u\i\ and C.P.~Plaut}
\address{ Valeri\u\i\ Nikolaevich Berestovski\u\i\ \\
Omsk Branch of Sobolev Institute of Mathematics SD RAS \\
644099, Omsk, ul. Pevtsova, 13, Russia}
\email{berestov@ofim.oscsbras.ru}
\address{Conrad Peck Plaut \\
Department of Mathematics\\
University of Tennessee\\
Knoxville, TN 37996, USA}
\email{cplaut@math.utk.edu}

\begin{abstract}
We prove that every length space $X$ is the orbit space (with the quotient
metric) of an $\mathbb{R}$-tree $\overline{X}$ via a free isometric action.
In fact, for many well-known spaces, such as connected complete Riemannian
manifolds of dimension at least two, the Menger sponge, and the Sierpin'ski
gasket and carpet, $\overline{X}$ is \textit{the} \textit{same}
\textquotedblleft universal\textquotedblright\ $\mathbb{R}$-tree $A_{%
\mathfrak{c}}$, which has valency $\mathfrak{c=}2^{\aleph _{0}}$ at each
point. The quotient mapping $\overline{\phi}:\overline{X}\rightarrow X$ is a
kind of generalized covering map called a URL-map, and $\overline{X}$ is the
unique (up to isometry) $\mathbb{R}$-tree that admits a URL-map onto $X$.
The map $\overline{\phi}$ is universal among URL-maps onto $X$ and in fact
is the \textquotedblleft mother of all metric universal
covers\textquotedblright\ in the following sense: All URL-maps, including
the traditional universal cover of a semi-locally simply connected length
space, may be naturally derived from it.

MSC Classification: 55Q05; 53C23, 28A80, 54F15
\end{abstract}

\keywords{universal cover, geodesic space, R-tree, fractal, continuum}
\maketitle

\section{Introduction and Main Results}

In this paper we construct the \textit{covering }$\mathbb{R}$\textit{-tree} $%
\overline{X}$ of a length space $X$, and the $\mathbb{R}$-tree universal
covering map $\overline{\phi }:\overline{X}\rightarrow X$. The function $%
\overline{\phi }$ is generally not a covering map in the traditional sense,
but shares important properties with metric covering maps. Recall that if $%
f:X\rightarrow Y$ is a traditional covering map and $Y$ is a length space,
the metric of $Y$ may be \textquotedblleft lifted\textquotedblright\ to $X$
in a unique way that makes $X$ a length space and $f$ a local isometry. The
function $f$ has two additional basic properties: (I) $f$ preserves the
length of rectifiable paths in the sense that $L(c)=L(f\circ c)$ for every
path $c$ in $X$ with finite length $L(c)$. (II) If $c$ is any rectifiable
path in $Y$ starting at a point $p$ and $f(q)=p$ then there is a unique path 
$c_{L}$ starting at $q$ such that $f\circ c_{L}=c$, and moreover $c_{L}$ is
rectifiable. A function $f$ between length spaces will be called \textit{%
unique rectifiable lifting (URL) }if $f$ has these two properties. Note that
a map between length spaces with condition (I) is known as an \textit{%
arcwise isometry} (\cite{Gr}). In fact any URL-map is a surjective open
arcwise isometry. The class of URL-maps is closed under composition and, as
we will see later, contains functions that are not local homeomorphisms--two
essential features if one has as a goal finding generalizations of the
traditional universal cover beyond spaces that are semi-locally simply
connected. In the case of a regular (or normal) covering map $f$, the deck
group $G$ of the covering map acts freely via isometries on $X$, and $Y$ is
the metric quotient $G\backslash X$. In other words, $Y$ is isometric to the
orbit space $G\backslash X$ with the Hausdorff metric on the orbits, and $f$
is the corresponding quotient map. In particular, the traditional universal
covering map $\phi :\widetilde{Y}\rightarrow Y$ (when it exists!) is regular
with deck group $\pi _{1}(Y)$. Recall that $\widetilde{Y}$ is unique (up to
isometry, with the lifted metric) and $\phi $ has a universal property: if $%
g:Z\rightarrow Y$ is a covering map then there is a unique (up to basepoint
choice) covering map $h:X\rightarrow Z$ such that $\phi =g\circ h$. For the
following theorem we define $\Lambda (X):=\lambda (X)/\eta (X)$, where $%
\lambda (X)$ is the group of rectifiable loops in the length space $X$
starting at a given basepoint and $\eta (X)$ is the normal subgroup of loops
that are homotopic \textit{in their image} to the trivial loop (see also
Definition \ref{lambdad} and Proposition \ref{hom}).

\begin{theorem}
\label{m1}For every length space $X$ there exists a unique (up to isometry) $%
\mathbb{R}$-tree $\overline{X}$, called the covering $\mathbb{R}$-tree, with
a URL-map $\overline{\phi }:\overline{X}\rightarrow X$. Moreover,

\begin{enumerate}
\item $\overline{X}$ is complete if and only if $X$ is complete.

\item If $Z$ is a length space and $f:Z\rightarrow X$ is a URL-map then
there is a unique (up to basepoint choice) URL-map $\overline{f}:\overline{X}%
\rightarrow Z$ such that $\overline{\phi }=f\circ \overline{f}$.

\item The group $\Lambda (X)$ acts freely via isometries on $\overline{X}$
with metric quotient map $\overline{\phi }:\overline{X}\rightarrow X=\Lambda
(X)\backslash\overline{X}$.
\end{enumerate}
\end{theorem}

The term \textquotedblleft $\mathbb{R}$-tree\textquotedblright\ was coined
by Morgan and Shalen (\cite{MS}) in 1984 to describe a type of space that
was first defined by Tits (\cite{T}) in 1977. Originally $\mathbb{R}$-trees
were defined as metric spaces with more than one point in which any two
points are joined by a unique geodesic, i.e. an arclength parameterized
curve with length equal to the distance between its endpoints. To avoid
trivial special cases, for this paper define \textquotedblleft length
space\textquotedblright\ (resp. \textquotedblleft geodesic
space\textquotedblright ) to be a metric space \textit{with at least two
points} such that any pair of points is joined by a path of length
arbitrarily close to the distance between them (resp. joined by a geodesic).
In the last three decades $\mathbb{R}$-trees have played a prominent role in
topology, geometry, and geometric group theory (see, for example, \cite{Be}, 
\cite{MS}). They are the most simple of geodesic spaces, and yet Theorem \ref%
{m1} shows that every length space, no matter how complex, is an orbit space
of an $\mathbb{R}$-tree.

In 1928 Menger asked whether (in modern terminology) every Peano continuum
(the continuous image of $[0,1]$) admits a compatible geodesic metric (\cite%
{M}). More than 20 years later the problem was given a positive answer
independently by Bing and Moise (\cite{B}, \cite{Moi}). In fact, these two
papers, together with earlier results of Menger, establish for compact,
connected metric spaces the equivalence of (1) local connectedness, (2)
local arcwise connectedness, and (3) the existence of a compatible geodesic
metric. A few years later, R. D. Anderson announced (\cite{An}) that every
Peano continuum is the continuous image of the Menger sponge $\mathbb{M}$
(the so-called \textquotedblleft universal curve\textquotedblright ) such
that each point pre-image is also $\mathbb{M}$. A proof of Anderson's
theorem was eventually published by Wilson in 1972 (\cite{W}), who at the
same time proved a strengthened conjecture of Anderson (\cite{An2}) by
showing that every Peano continuum is the image of a $\mathbb{M}$ via a
light open mapping (\textit{light }means every point pre-image is totally
disconnected). Anderson's conjecture was part of a long-standing effort to
construct dimension-raising open mappings, beginning with an example of
Kolmogorov in 1937 (\cite{K}) from a $1$-dimensional Peano continuum to a $2$%
-dimensional space (see also \cite{Ke} for a dimension-raising light open
mapping). Proposition \ref{lite} and the fact that a (non-trivial) $\mathbb{R%
}$-tree $X$ is simply connected with small inductive dimension $ind(X)=1$ (%
\cite{AB}) give us:

\begin{corollary}
\label{pc}Every non-trivial space admitting a compatible length metric is
the image via a light open mapping of a simply connected space $X$ with $%
ind(X)=1$.
\end{corollary}

Consider the following fractals: the Sierpin'ski carpet $S_{c}$, the
Sierpin'ski gasket $S_{g}$, or $\mathbb{M}$. As is well-known, each such
space $X$ admits a geodesic metric $d$ bi-Lipshitz equivalent to the metric
induced by the metric $\rho $ of the ambient Euclidean space; $d(x,y)$ is
defined as the infimum of the length of paths in $X$ joining the points $x$
and $y$, where the length is measured in the metric $\rho $. An $\mathbb{R}$%
-tree that appears very naturally in our work is the $\mathfrak{c}$%
-universal $\mathbb{R}$-tree $A_{\mathfrak{c}}$ introduced in \cite{MNO},
which has valency $\mathfrak{c}=2^{\aleph _{0}}$ (cardinality of the
continuum) at each point. $A_{\mathfrak{c}}$ was shown in \cite{MNO} to be
metrically homogeneous, and \textquotedblleft universal\textquotedblright\
in the sense that every $\mathbb{R}$-tree of valency at most $\mathfrak{c}$
isometrically embeds in $A_{\mathfrak{c}}$. This is analogous to the
original way in which $\mathbb{M}$ was considered \textquotedblleft
universal\textquotedblright .

\begin{theorem}
\label{universal} If $X$ is a separable length space, then $\overline{X}$ is
a sub-tree of $A_{\mathfrak{c}}$. If in addition $X$ is complete and
contains a bi-Lipschitz copy of $S_{g}$ or $S_c$ at every point, e.g. if $X$
is $S_{c} $, $S_{g}$, $\mathbb{M}$, or a complete Riemannian manifold of
dimension at least two, then $\overline{X}$ is isometric to $A_{\mathfrak{c}%
} $.
\end{theorem}

Put another way, every separable length space may be obtained by starting
with a subtree of $A_{\mathfrak{c}}$ and taking a quotient of that subtree
via a free isometric action. Another consequence of this theorem is an
explicit construction of $A_{\mathfrak{c}}$ starting with any of the above
spaces (see the proof of Theorem \ref{m1}). Our results, combined with the
Anderson-Wilson Theorem, show that $A_{\mathfrak{c}}$ is \textquotedblleft
universal\textquotedblright\ in another sense, which is similar to the
second way in which $\mathbb{M}$ may be regarded as \textquotedblleft
universal\textquotedblright :

\begin{corollary}
Every Peano continuum is the image of $A_{\mathfrak{c}}$ via a light open
mapping.
\end{corollary}

Finally, Theorem \ref{m1} provides yet a third (categorial) way in which $A_{%
\mathfrak{c}}$ may be considered \textquotedblleft
universal\textquotedblright .

In addition to the properties described in the first paragraph, the
traditional universal covering has two useful properties: the universal
covering map is a fibration and the universal covering space is simply
connected. These properties are used to classify traditional covering maps
according to the subgroups of the fundamental group whose (representative)
elements lift as loops. While $\overline{X}$ is simply connected, $\overline{%
\phi }$ is generally not a fibration. In fact, the only homotopies that may
be lifted must already be \textquotedblleft tree-like\textquotedblright ,
and it would be interesting to understand which homotopies do lift. For
example, Piotr Hajlasz and Jeremy Tyson have constructed Lipschitz
surjections from cubes onto compact, quasiconvex, doubling metric spaces (%
\cite{HT}), and it is easy to see from their construction that these
mappings lift to the covering $\mathbb{R}$-tree. The mapping $\overline{\phi 
}$ is not only \textit{not} a fibration, it is in some sense as far as
possible from being a fibration. This turns out to be an advantage. For any
subgroup $G$ of $\Lambda (Y)$ we define a kind of \textquotedblleft designer
fibration\textquotedblright\ called a URL($G$)-map $f:X\rightarrow Y$, which
is a URL-map that lifts any (representative) loop in $G$ as a loop. If the
loops in $G$ are the \textit{only }loops that lift as loops then $f$ is
called \textquotedblleft $G$-universal\textquotedblright\ (Definition \ref%
{gsimp}). We prove:

\begin{theorem}
\label{m2}Let $X$ be a length space and $G$ be a subgroup of $\Lambda (X)$
considered as a group of isometries of $\overline{X}$. Suppose that the
orbits of $G$ are closed, the orbit space $G\backslash \overline{X}$ is
supplied with the quotient metric, and $\psi _{G}:\overline{X}\rightarrow
G\backslash \overline{X}:=\overline{X}^{G}$ is the quotient mapping. Then
there exists unique (continuous) map $\overline{\phi }^{G}:\overline{X}%
^{G}\rightarrow X$ such that $\overline{\phi }^{G}\circ \psi _{G}=\overline{%
\phi }$. Moreover, if $\overline{\phi }^{G}$ is a URL-map then

\begin{enumerate}
\item $\overline{X}^{G}$ is the unique (up to isomorphism) $G$-universal
space with a URL($G$)-map $\overline{\phi }^{G}:\overline{X}^{G}\rightarrow
X $ (also unique up to basepoint choice).

\item For any URL($G$)-map $f:Y\rightarrow X$ there is a unique (up to
basepoint choice) URL-map $f^{G}:\overline{X}^{G}\rightarrow Y$ such that $%
\overline{\phi }^{G}=f\circ f^{G}$.

\item If $G$ is a normal subgroup of $\Lambda (X)$ then the group $\pi
_{1}^{G}(X):=\Lambda (X)/G$ acts via isometries on $\overline{X}^{G}$ and $%
\overline{\phi }^{G}$ is the quotient mapping with respect to this action.

\item $X$ is complete if and only if $\overline{X}^{G}$ is complete.
\end{enumerate}
\end{theorem}

When $G$ is normal we call $\pi _{1}^{G}(X)$ the $G$-fundamental group of $X$%
. The next corollary shows that $\overline{\phi }^{G}$ is a true
generalization of the traditional metric universal covering.

\begin{corollary}
\label{slsc}Let $X$ be a semilocally simply connected length space. Then $%
\overline{\phi }^{H_{T}}$ is the traditional metric universal covering of $X$
and $\pi _{1}^{H_{T}}(X)$ is naturally isomorphic to $\pi _{1}(X)$.
\end{corollary}

In general, there is a natural homomorphism $h_{\Lambda }:\Lambda
(X)\rightarrow \pi _{1}(X)$, the kernel of which is $H_{T}$ and the image of
which is the subgroup $\mu _{1}(X)$ of $\pi _{1}(X)$ consisting of those
homotopy classes having a rectifiable representative (see Proposition \ref%
{fund}). Length spaces with bad (or unknown) local topology appear often in
geometry and topology, from the classical examples such as the Menger sponge
and Hawiian earring to Gromov-Hausdorff limits of Riemannian manifolds. Very
recently, Sormani and Wei have shown that such limits of manifolds with
non-negative Ricci curvature have a universal covering in the categorical
sense, but it is unknown whether this covering is simply connected (\cite%
{SWHC}, \cite{SWNC}). Their work, in turn, was partly motivated by the
40-year old conjecture of Milnor that such non-positively curved manifolds
have finitely generated fundamental groups (\cite{Mi}). These papers and
others (\cite{Wy}, \cite{SWCC}, \cite{SWCS}) involve studying covering maps
determined by geometrically significant groups of loops through a
construction of Spanier (\cite{Sp}). Our paper provides a much more general
framework for such efforts. In fact, for any open cover $\mathcal{U}$ of a
connected, locally arcwise connected space, Spanier defined a covering space
corresponding to a certain subgroup $G_{\mathcal{U}}$ of the fundamental
group. When $X$ is a length space, Spanier's covering map is of the form $%
\overline{\phi }^{G}$ where $G:=h_{\Lambda }^{-1}(G_{\mathcal{U}})$.

In \cite{BPU} we introduced the uniform universal covering (UU-covering) $%
\phi :\widetilde{X}\rightarrow X$. The UU-covering is an analog of the
universal covering for a large class of uniform spaces called \textit{%
coverable} spaces, which includes all spaces admitting a length metric, and
therefore all Peano continua. However, the UU-covering is not always
satisfactory from a geometric standpoint. For example, the UU-coverings of
the Hawaiian earring and the Menger sponge are connected but not arcwise
connected, and the UU-covering of a non-compact semi-locally simply
connected length space may not be its traditional universal covering (see 
\cite{BPU}). In contrast, the covering $\mathbb{R}$-tree involves the more
specialized (but extremely important) class of length spaces, but provides a
way to construct a variety of generalized universal metric covering spaces
that are always length spaces. Even so, the UU-covering proves useful in a
couple of ways in the current paper. We use it to prove two essential
theorems about paths (Theorems \ref{normal} and \ref{wnn}), and we use it to
obtain additional examples of URL-maps that are not traditional covering
maps. When $X$ is a length space we show there exists a \textquotedblleft
metric core\textquotedblright\ $\widehat{X}\subset \widetilde{X}$ that is a
length space. Moreover, the restriction $\widehat{\phi }:\widehat{X}%
\rightarrow X$ of the UU-covering map is a URL-map and a \textquotedblleft
metric fibration\textquotedblright\ in the sense that Lipschitz maps from
geometrically reasonable simply connected domains may be lifted (Theorem \ref%
{core}). This may prove useful since in the case of the above-mentioned
fractals and the compact separable infinite dimensional torus $T^{\infty }$,
the metric core is a $CAT(0)$ space even though the original space is not
even locally simply connected. In fact the metric core of $T^{\infty }$ is
separable Hilbert space. In the case of a uniformly $1$-dimensional length
space, $\overline{X}$ is naturally isometric to $\widehat{X}$ (Corollary \ref%
{same}). For length spaces (as opposed to uniform spaces in general), the
construction (if not the proofs!) of the UU-covering is simpler to describe
and we may use simplified notation. We give such a description and establish
our notation in a short appendix.

\section{The covering $\mathbb{R}$-tree\label{p1sec}}

The following are equivalent for a geodesic space $X$ (see \cite{Mo}, \cite%
{BH}, \cite{AB}): (1) $X$ is an $\mathbb{R}$-tree. (2) $X$ is $0$-hyperbolic
in Gromov's sense. (3) $X$ contains more then one point and is $CAT(K)$%
-space for all $K\leq 0$. (4) $X$ is simply connected and its
small-inductive Urysohn-Menger dimension is 1. See, for example, \cite{BH}
for the definitions of $0$-hyperbolic and $CAT(K)$-space; we will not need
the definitions in this paper. However, note that a corollary of Theorem \ref%
{m1} is that every length space is the metric quotient of a $CAT(K)$
geodesic space for any $K\leq 0$.

A few words about dimension are in order. For a metric space $X$ we denote
the small (resp. large) inductive dimension by $ind(X)$ (resp. $Ind(X)$),
and covering dimension by $\dim (X)$. It is known that for an arbitrary
metric space $X$, there are the Katetov equality $Ind(X)=\dim (X)$ and the
inequality $ind(X)\leq Ind(X)$, see \cite{AP}. If $X$ is also separable (in
particular compact) then $ind(X)=Ind(X)=dim(X)$, \cite{HW}. Generally an $%
\mathbb{R}$-tree is not separable (for example, $A_{\mathfrak{c}}$). For
metric spaces (or more generally uniform spaces), one may also consider
various definitions of \textquotedblleft uniform dimension\textquotedblright
. We use the same definition here as in \cite{BPU}, which is the same as the
definition of covering dimension for a topological space except that the
open covers involved are \textquotedblleft uniform\textquotedblright\ in the
following sense: An open cover of a metric space is uniform if it is refined
by a cover of $\varepsilon $-balls for some $\varepsilon >0$. This dimension
is called \textquotedblleft large dimension\textquotedblright\ in \cite{I}
and we will denote it by $u\dim (X)$ in this paper. In general we do not
know precisely the relationship between uniform dimension and the various
topological notions of dimension.

\begin{notation}
We will write \textquotedblleft $X$ has dimension $\leq n$ in some sense
(resp. in every sense)\textquotedblright\ if $ind(X)$, $Ind(X)$, $\dim (X)$,
or (resp. and) $u\dim (X)$ is at most $n$.
\end{notation}

For example, it is known that if $X$ has dimension $\leq n$ in some sense
and $C\subset X$ is compact then $C$ has dimension $\leq n$ in every sense.
We will use this fact frequently below.

We next recall some basic background about paths, by which we mean
continuous functions from compact intervals into metric spaces. While usage
varies in the literature, we will say that two paths $c_{1},c_{2}:I%
\rightarrow X$ are Fr\'{e}chet equivalent if there exists an
order-preserving homeomorphism $c:I\rightarrow I$ such that $%
c_{1}=c_{2}\circ c$. If $c$ is a path in a metric space we denote by $L(c)$
its length, writing $L(c)=\infty $ if $c$ is not rectifiable. An arclength
parameterized path is $1$-Lipschitz, hence does not increase Hausdorff
dimension (which is not smaller than covering dimension). It follows that
the image of any nonconstant rectifiable path is a uniformly $1$-dimensional
Peano continuum. We will consider parameterizations only up to Fr\'{e}chet
equivalence, which simplifies many discussions. For example, given paths $%
c_{1}$ and $c_{2}$ defined on $[0,L_{1}]$ and $[0,L_{2}]$, respectively,
such that the starting point of $c_{2}$ is the endpoint of $c_{1}$, we will
simply refer to the concatenation $c_{1}\ast c_{2}$ on $[0,L_{1}+L_{2}]$
without mentioning the linear reparameterization of $c_{2}$ to $%
[L_{1},L_{1}+L_{2}]$ that is technically required to concatenate them. For
any path $c$ we will denote by $c^{-1}$ the same path with orientation
reversed (again using any convenient parameterization in the Fr\'{e}chet
equivalence class). When it comes to homotopies, we may say that two paths $%
c_{1}$ and $c_{2}$ are \textquotedblleft
endpoint-homotopic\textquotedblright\ even though, strictly speaking, there
is only a homotopy between paths $c_{1}^{\prime }$ and $c_{2}^{\prime }$
that are Fr\'{e}chet equivalent to $c_{1}$ and $c_{2}$, respectively, that
share a common domain.

Recall that in \cite{CF} a path $c:[a,b]\rightarrow X$ in a metric space $X$
is called \textit{normal} if there is no nontrivial subsegment $%
J=[u,v]\subset I$ such that $c(u)=c(v),$ and $c|J$ is homotopic to a
constant relative to $\{u,v\}$.

\begin{definition}
\label{cfd1}A path $c:[a,b]\rightarrow X$ in a metric space $X$ is called
weakly normal if there is no nontrivial subsegment $J=[u,v]\subset \lbrack
a,b]$ such that $c(u)=c(v)$ and $c|_{J}$ is homotopic in $c(J)$ to a
constant relative to $\{u,v\}$. A weakly normal (rectifiable) arclength
parameterized path (resp. loop) $c:[a,a+L]\rightarrow X$ in a metric space $%
X $ will be called a $\rho $-path (resp. $\rho $-loop).
\end{definition}

Put another way, a weakly normal path is one that is normal in its own
image; clearly every normal path is weakly normal.

\begin{remark}
In order to avoid special cases below, we permit the domain of a constant
function to be an interval of the form $[a,a]=\{a\}$; such a
parameterization of a constant map is clearly a $\rho $-path.
\end{remark}

A version of the next theorem was proved for $1$-dimensional separable
metric spaces using very different methods (including Zorn's Lemma) in \cite%
{CF} (Lemma 3.1 and Theorem 3.1), also without considering lengths of paths.
For our purposes it is very important to consider lengths and
rectifiability. Moreover, as we have pointed out earlier, the $\mathbb{R}$%
-trees we will be considering may not be separable and therefore it is
important to remove this assumption.

\begin{theorem}
\label{normal}Each fixed-endpoint homotopy class of a path $c:I\rightarrow X$
in metric space $X$ that is $1$-dimensional in some sense contains a unique
(up to Fr\'{e}chet equivalence) normal path $c_{n}$. In addition, $%
L(c_{n})\leq L(c)$.
\end{theorem}

\begin{proof}
First suppose that $X$ is a Peano continuum. Then $X$ has a UU-covering $%
\widetilde{X}$. Choose as a basepoint $\ast $ the starting point of $c$ and
suppose first that $c$ is not a null-homotopic loop. By the lifting property
of the UU-covering (see the Appendix), there is unique path $\tilde{c}%
:I\rightarrow \widetilde{X}$ starting at $\ast $ such that $\phi \circ 
\tilde{c}=c$. Since $X$ is compact, $X$ is uniformly $1$-dimensional and we
may apply Proposition \ref{inj} to conclude that the function $\lambda $
described in the appendix is injective. This means that $\tilde{c}$ is a
loop if and only if $c$ is a null-homotopic loop. Therefore by our initial
assumption $\tilde{c}$ is not a loop. The image $C=\tilde{c}(I)$ is a Peano
continuum that contains no topological circle by Proposition \ref{inj}. By
the Hahn-Mazurkiewicz-Sierpin'ski Theorem, $C$ is locally connected, hence a
dendrite, see Section 51, VI of \cite{Ku}. Then $C$ is contractible by
Corollary 7 in Section 54, VII of \cite{Ku}. Hence there is a fixed-endpoint
homotopy $h:I\times I\rightarrow C$ such that $h(\cdot ,0)=\tilde{c}$, and $%
h(\cdot ,1)$ is a topological embedding $\tilde{c}_{1}$, whose image is the
unique arc $a$ in $C$ joining $\ast \in C$ and $\tilde{c}(1)=\tilde{c}%
_{1}(1) $ (this arc exists by Corollary 7, Section 51, VI in \cite{Ku}). By
Proposition \ref{inj}, $a$ is also the unique arc in $\widetilde{X}$ joining 
$\widetilde{c}(1)$ and $\ast $. Then $\phi \circ h$ is a fixed-endpoint
homotopy in $C$ from $c$ to the unique (up to Fr\'{e}chet equivalence)
normal path $c_{n}:=\phi \circ \tilde{c}_{1}$. To finish the proof in the
case when $c$ is not a null-homotopic loop, we need to prove that $%
L(c_{n})\leq L(c)$ if $c$ is rectifiable. By the previous argument, $%
a\subset C=\tilde{c}(I)$, hence $c_{n}(I)\subset c(I)$. Let $%
t_{0}=0<t_{1}<\dots <t_{m}=1$ be a partition of $I$ and $%
x_{0}=x=c_{n}(t_{0}),x_{1}=c_{n}(t_{1}),\dots ,x_{m}=c_{n}(t_{m})$. Then for
any $i=1,\dots ,m$ there is a maximal $s_{i}\in I$ such that $%
x_{i}=c(s_{i}). $ It is clear that $s_{0}=0<s_{1}<\dots <s_{m}=1$, so we get
another partition of $I$. Now we have $%
\sum_{i=0}^{n-1}d(c_{n}(t_{i}),c_{n}(t_{i+1}))=%
\sum_{i=0}^{n-1}d(c(s_{i}),c(s_{i+1}))\leq L(c)$, and $L(c_{n})\leq L(c)$.

If $c$ is a non-constant null-homotopic loop then we may write it as the
concatenation of two paths that are not null-homotopic loops and apply the
case just considered; the path $c_{n}$ in this case is constant. If $c$ is
constant, its homotopy class could not contain a non-constant normal path by
the preceding sentence.

Returning to the general case for $X$, note that the image $C$ of $c$ is a
Peano continuum that is $1$-dimensional (with the subspace metric) in the
same sense that $X$ is, and hence we may apply the above special case to
obtain the existence of the desired path $c_{n}$ in the fixed-endpoint
homotopy class of $c$ in $C$. Now suppose that $c^{\prime }$ is any normal
path in the fixed-endpoint homotopy class of $c$ in $X$. Let $h$ be a
homotopy from $c_{n}$ to $c^{\prime }$. We may apply a similar argument
using the compact image $K$ of $h$ to conclude from the uniqueness in the
case of a Peano continuum that $c_{n}$ and $c^{\prime }$ are Fr\'{e}chet
equivalent.
\end{proof}

\begin{theorem}
\label{wnn} A path $c:I\rightarrow X$ in a metric space that is $1$%
-dimensional in some sense is normal if and only if $c$ is weakly normal.
\end{theorem}

\begin{proof}
As in the proof of Theorem \ref{normal}, we may assume that $X$ is a Peano
continuum. We have already observed that every normal path is weakly normal.
Let us suppose that the path $c$ is not normal. Then there is a nontrivial
subsegment $J=[u,v]\subset I$ such that $c(u)=c(v),$ and $c|J$ is homotopic
to a constant ($c(u)=c(v)$) relative to $\{u,v\}$. Let $h:J\times
I\rightarrow X$ be the corresponding homotopy. By the lifting property,
there is a unique homotopy $\tilde{h}:J\times I\rightarrow \tilde{X}$ such
that $\tilde{h}(u,0)=\tilde{c}(u)$ and $h=\phi \circ \tilde{h}.$ Then $\phi
\circ \tilde{h}(u,s)=\phi \circ \tilde{h}(v,s)=c(u)=c(v)$ for all $s\in I,$
and by Remark \ref{light}, $\tilde{h}(u,s)=\tilde{c}(u)$ and $\tilde{h}(v,s)=%
\tilde{c}(v)$ for all $s\in I$. Evidently, $\tilde{h}(t,1)$ is constant for
all $t\in J,$ so $\tilde{c}(u)=\tilde{c}(v).$ Then the restriction $\tilde{c}%
|J$ is a loop in $\tilde{X}$ and its image $C=\tilde{c}(J)$ is a dendrite as
in the proof of Theorem \ref{normal}. There is a fixed-endpoint homotopy $%
\tilde{h}_{1}:J\times I\rightarrow C$ such that $\tilde{h}(\cdot ,0)=\tilde{c%
}|J$ and $\tilde{h}_{1}(t,1)\equiv \tilde{c}(u)=\tilde{c}(v)$ for all $t\in
J $. Then $h_{1}=\phi \circ \tilde{h}_{1}$ gives a fixed-endpoint homotopy
in $c(J)$ to the constant path $c(u)=c(v)$. This means that the path $c$ is
not weakly normal.
\end{proof}

From our previous observation that the image of a rectifiable path is a $1$%
-dimensional Peano continuum we may conclude:

\begin{corollary}
\label{rep}Let $X$ be a metric space and $c$ be path in $X$, with image $C$
that is $1$-dimensional in some sense (resp. rectifiable). There is a unique
weakly normal (resp. $\rho $-path) $\gamma :[0,L(\gamma )]\rightarrow C$
that is fixed-endpoint homotopic to $c$ in $C$. In addition, $L(c)\geq
L(\gamma )$.
\end{corollary}

\begin{corollary}
In any length space, for every pair of points $p,q$, $d(p,q)$ is the infimum
of the lengths of injective $\rho $-paths joining $p$ and $q$.
\end{corollary}

Note that geodesics are injective and hence normal.

\begin{corollary}
\label{runiq}If $X$ is a simply connected length space that is $1$%
-dimensional in some sense (in particular if $X$ is an $\mathbb{R}$-tree)
then every pair of points is joined by a unique $\rho $-path, which is also
a geodesic (hence the unique geodesic) joining them. Moreover, this $\rho $%
-path is contained in the image of every path joining the two points.
\end{corollary}

\begin{proof}
There is such a path, since any two points $p,q$ are connected by a $\rho $%
-path. Moreover, since $X$ is simply connected, any two $\rho $-paths
joining them must be fixed-endpoint homotopic, hence by Theorems \ref{normal}
and \ref{wnn}, there can only be one $\rho $-path joining them. Finally, $%
d(p,q)$ is equal to the infimum of the lengths of $\rho $-paths joining
them--but up to Fr\'{e}chet equivalence there is only one such path, so the
infimum must be a minimum.
\end{proof}

Note that the concatenation of a $\rho $-path $c$ followed by a $\rho $-path 
$d$ may not be a $\rho $-path To resolve this problem we define the
\textquotedblleft cancelled concatenation\textquotedblright\ $c\star d$ to
be the unique $\rho $-path in the fixed-endpoint homotopy class of the
concatenation $c\ast d$, in the image of $c\ast d$ (Corollary \ref{rep}). By
uniqueness, we see more concretely that $c\star d$ is obtained from $c\ast d$
by removing the maximal final segment of $c$ that coincides with an initial
segment of $d$ with reversed orientation, and removing that initial segment
of $d$ as well. Also by uniqueness, the associative law $(a\star b)\star
c=a\star (b\star c)$ is satisfied. Moreover, cancelled concatenation on $%
\rho $-loops at a fixed basepoint is a group operation, where the constant
loop is the identity and the inverse of $c$ is the $\rho $-loop $c^{-1}$
already discussed above.

\begin{definition}
\label{lambdad}Let $X$ be a length space with basepoint $\ast $ and define $%
\left( \overline{X},\ast \right) $ (or simply $\overline{X}$ when no
confusion with regard to basepoint will result) to be the set of all $\rho $%
-paths $c:[0,L(c)]\rightarrow X$ starting at $\ast $ (i.e. $c(0)=\ast $),
and let $\overline{\phi }:\overline{X}\rightarrow X$ be the endpoint
mapping. We place the following metric on $\overline{X}$: Let $c_{1}\wedge
c_{2}:[0,b]\rightarrow X$ be the restriction of $c_{1}$ (and $c_{2}$) to the
largest interval $[0,b]$ on which $c_{1}$ and $c_{2}$ coincide. For $%
c_{1},c_{2}\in \overline{X}$, define 
\begin{equation}
d_{\overline{X}}(c_{1},c_{2}):=L(c_{1})+L(c_{2})-2L(c_{1}\wedge
c_{2})=L(c_{1}^{-1}\star c_{2})  \label{d}
\end{equation}%
We let $\Lambda (X,\ast )\subset \left( \overline{X},\ast \right) $ (or
simply $\Lambda (X)$) denote the group of all $\rho $-loops starting at $%
\ast $, with the subspace metric.
\end{definition}

It is easy to check that $d_{\overline{X}}$ is a metric.

\begin{proposition}
\label{basepoint}For a length space $X$ with basepoints $\ast $ and $\ast
^{\prime }$, any $\rho $-path $k$ from $\ast ^{\prime }$ to $\ast $ induces
an isometry from $(\overline{X},\ast )$ to $(\overline{X},\ast ^{\prime })$
defined by $c\mapsto k\star c$, and composition of this isometry with the
endpoint mapping coincides with $\overline{\phi }$.
\end{proposition}

\begin{proof}
In fact, for $\rho $-curves $c_{1},c_{2}$ starting at $\ast $, the
previously discussed associative law implies 
\begin{equation*}
d_{\overline{X}}(k\star c_{1},k\star c_{2})=L(\left( k\star
c_{1}\right)^{-1} \star \left( k\star c_{2}\right))=L(c_{1}^{-1}\star
c_{2})=d_{\overline{X}}(c_{1},c_{2})
\end{equation*}
\end{proof}

That is, up to isometry, $\overline{\phi }:\overline{X}\rightarrow X $ is
independent of choice of basepoint in $X$. For this reason we will often
avoid discussion of basepoints and simply assume that all functions are
basepoint preserving, in particular choosing the constant path $\ast $ as
the basepoint in $\overline{X}$.

\begin{corollary}
\label{imp} The group $\Lambda (X)$ acts freely as isometries on $(\overline{%
X},\ast )$ (hence on $\Lambda (X)$) by left multiplication.
\end{corollary}

\begin{lemma}
For any length space $X$, the metric on $\Lambda (X)$ is complete, hence $%
\Lambda (X)$ is a complete topological group.
\end{lemma}

\begin{proof}
First note that it is an immediate consequence of the definition of the
metric that for any $c,d\in \overline{X}$, 
\begin{equation*}
\left\vert L(c)-L(d)\right\vert \leq d(c,d)
\end{equation*}%
Let $\left( c_{i}\right) $ be a Cauchy sequence in $\Lambda (X)$ with $c_{i}$
parameterized on $[0,L_{i}]$. By the above inequality, $(L_{i})$ is a real
sequence converging to a non-negative real number $L$. Now for any $%
L^{\prime }<L$ the loops $c_{i}$ are defined on $[0,L^{\prime }]$ and equal,
for all sufficiently large $i$. We may define $c:[0,L]\rightarrow X$ by $%
c(t)=c_{i}(t)$ for all large $i$, and $c(L)=\ast $. It is easy to check that
this defines an element of $\Lambda (X)$, and certainly $c_{i}\rightarrow c$
in the metric of $\Lambda (X)$. It is a standard result that a group with a
complete left-invariant metric is complete as a topological group.
\end{proof}

\begin{remark}
\label{basepoint1}It should be observed that the isometry from Proposition %
\ref{basepoint} does not take loops to loops and so does not induce an
isometry from $\Lambda (X,\ast )$ to $\Lambda (X,\ast ^{\prime })$. In fact,
these groups are in general isomorphic as groups, but not isometric. An
isomorphism from $\Lambda (X,\ast )$ to $\Lambda (X,\ast ^{\prime })$ can be
obtained by choosing any $\rho $-path $c$ from $\ast ^{\prime }$ to $\ast $
and taking $d\in \Lambda (X,\ast )$ to $c\star d\star c^{-1}\in \Lambda
(X,\ast ^{\prime })$. To see that these groups need not be isometric, let $X$
consist of a circle of circumference $1$ having a segment $[a,b]$ of length $%
1$ glued to it at $a$. The shortest distance from a non-trivial element of $%
\Lambda (X,a)$ to the trivial loop at $a$ is $1$, but the shortest distance
from a non-trivial element of $\Lambda (X,b)$ to the trivial loop at $b$ is $%
2$.
\end{remark}

One can easily prove the following:

\begin{proposition}
\label{hom}Let $\lambda (X)$ be the collection of all rectifiable loops at
the basepoint in $X$ with concatenation as the group operation. Define $%
h:\lambda (X)\rightarrow \Lambda (X)$ by letting $h(c)$ be the unique $\rho $%
-path in the fixed-endpoint homotopy class of $c$ in the image of $c$. Then $%
h$ is a surjective homomorphism with kernel equal to the (normal) subgroup $%
\eta (X)$ of all null-homotopic in itself rectifiable loops. That is, $%
\Lambda (X)=\lambda (X)/\eta (X)$.
\end{proposition}

For the next theorem we need to recall some definitions (see the survey
article \cite{Po} for more details). Submetries were introduced by the first
author as a generalization of the notion of Riemannian submersion (see \cite%
{BG} and references there). Recall that if $X,Y$ are metric spaces, $%
f:X\rightarrow Y$ is a submetry (resp. weak submetry) if for every closed
ball $B(p,r)$ (resp. open ball $U(p,r)$) in $X$, $f(B(p,r))=B(f(p),r)$
(resp. $f(U(p,r))=U(f(p),r)$). There are two trivially equivalent conditions
that we will use without further comment: (1) The function $f$ is distance
non-increasing (i.e. $1$-Lipschitz) and for every $x,y\in Y$ and $z\in
f^{-1}(x)$ (resp. and any $\varepsilon >0$) there exists $w\in f^{-1}(y)$
such that $d(z,w)=d(x,y)$ (resp. $d(z,w)<d(x,y)+\varepsilon $). (2) For
every $x,y\in X$ and $z\in f^{-1}(x)$, $d(x,y)=\min \{d(z,w):w\in f^{-1}(y)\}
$ (resp. $d(x,y)=\inf \{d(z,w):z\in f^{-1}(x)$ and $w\in f^{-1}(y)\}$. An
easy consequence of the definition is that $f$ is open. Obviously weak
submetries do not increase the lengths of paths. It follows that if $%
f:X\rightarrow Y$ is a weak submetry and $X$ is a length space then so is $Y$%
, although $Y$ need not be a geodesic space even if $X$ is. Every isometry
is obviously a submetry; on the other hand any injective weak submetry is an
isometry. If a group $G$ acts on a metric space $M$ by isometries then the
orbit space $G\backslash M$ may be given the quotient pseudometric, namely
the Hausdorff distance between the orbits (which may be $0$ for different
orbits if the two orbits are not closed). In fact, $G\backslash M$ is a
metric space if and only if the orbit $Gx=\{g(x):g\in G\}$ of every $x\in M$
is closed. Since $G$ acts by isometries and is transitive on each orbit, we
also have $d(Gx,Gy)=\inf \{d(g(x),y):g\in G\}$. It follows from the second
characterization above that with respect to this metric the \textit{metric
quotient mapping} $\phi :M\rightarrow G\backslash M$ is a weak submetry.

\begin{theorem}
\label{ut}For any length space $X$, $\overline{X}$ is an $\mathbb{R}$-tree
and $\overline{\phi }$ is the metric quotient mapping (hence a weak
submetry) with respect to the isometric action of $\Lambda (X)$.
\end{theorem}

\begin{proof}
We will show that $\overline{X}$ is an $\mathbb{R}$-tree using the
characterization (2) from the beginning of this section. We will also denote
by $\ast $ the element of $\overline{X}$ that is simply the constant path at 
$\ast \in X$. Let $c_{1},c_{2}\in \overline{X}$, defined on $[0,L_{1}]$, $%
[0,L_{2}]$, respectively. Let 
\begin{equation*}
s_{0}:=\max \{s:c_{1}(t)=c_{2}(t)\text{ for all }t\in \lbrack 0,s]\}
\end{equation*}%
and define $C(s)$ for $s\in \lbrack 0,L_{1}+L_{2}-2s_{0}]$ as follows. For $%
s\in \lbrack 0,L_{1}-s_{0}]$ let $C(s)$ be the restriction of $c_{1}$ to $%
[0,L_{1}-s]$. For $s\in \lbrack L_{1}-s_{0},L_{1}+L_{2}-2s_{0}]$ let $C(s)$
be the restriction of $c_{2}$ to $[0,s-L_{1}+2s_{0}]$. Certainly $C(s)$ is a
geodesic in $\overline{X},$ joining $c_{1}$ and $c_{2}$. This implies that $%
\overline{X}$ is a geodesic space.

We see from Formula (\ref{d}) that the so-called Gromov product 
\begin{equation*}
(c_{1},c_{2})_{\ast }:=\frac{1}{2}[d_{\overline{X}}(\ast ,c_{1})+d_{%
\overline{X}}(\ast ,c_{2})-d_{\overline{X}}(c_{1},c_{2})]
\end{equation*}%
with respect to the point $\ast $ (see, for example, \cite{BH}, page 410) is
equal to $L(c_{1}\wedge c_{2})$. Also we see immediately that $c_{1}\wedge
c_{2}$ contains as a subpath $(c_{1}\wedge c_{3})\wedge (c_{2}\wedge c_{3})$
for any $c_{3}$. Then it follows from these two statements that 
\begin{equation}
(c_{1},c_{2})_{\ast }\geq \min \{(c_{1},c_{3})_{\ast },(c_{3},c_{2})_{\ast
}\}  \label{c}
\end{equation}%
for any $c_{1},c_{2},c_{3}\in \overline{X}$. This means that $\overline{X}$
is 0-hyperbolic \textquotedblleft with respect to the point $\ast $%
\textquotedblright , whereas 0-hyperbolicity itself means that the equation (%
\ref{c}) must be always satisfied if we replace $\ast $ by any path-point $c$%
. By Remark 1.21 (page 410 of \cite{BH}), $\overline{X}$ is actually
0-hyperbolic, so it is an $\mathbb{R}$-tree.

To show that $\overline{\phi }$ is the quotient mapping we need to check:
(1) For any $x\in X$, $\overline{\phi }^{-1}(x)$ is precisely the orbit $%
\Lambda (X)y$ of any point $y\in \overline{\phi }^{-1}(x)$ and (2) for any $%
x,w\in X$, $d(x,w)$ is equal to the Hausdorff distance $d_{H}(\overline{\phi 
}^{-1}(x),\overline{\phi }^{-1}(w))$. Note that if $\overline{\phi }(c)=%
\overline{\phi }(d)$ then by definition $c$ and $d$ have the same endpoint.
So $c\star d^{-1}\in \Lambda (X)$ and $c=(c\star d^{-1})\star d$, which
means that $c$ and $d$ lie in the same orbit. Conversely, if $c$ and $d$ lie
in the same orbit, $c=k\star d$ for some $k\in \Lambda (X)$, so both curves
have the same endpoint. This proves condition (1). Now by definition, for
any $c,f\in \overline{X}$, $d(c,f)$ is at least the distance between their
endpoints $p,q$, respectively, since $d(c,f)$ is the length of a curve
joining $p$ and $q$. This implies that $d_{H}(\overline{\phi }^{-1}(p),%
\overline{\phi }^{-1}(q))\geq d(p,q)$. On the other hand, given any $p,q\in
X $ and $\varepsilon >0$ there is a $\rho $-path $c$ from $\ast $ to $p$ and
a $\rho $-path $k$ from $p$ to $q$ such that $L(k)\leq d(p,q)+\varepsilon $.
Let $d:=c\star k$. Then $d(c,f)=L(f^{-1}\star c)\leq L(k)\leq
d(p,q)+\varepsilon $. Letting $\varepsilon \rightarrow 0$ we see $d_{H}(%
\overline{\phi }^{-1}(p),\overline{\phi }^{-1}(q))\leq d(p,q)$.
\end{proof}

\begin{proposition}
\label{complete}If $f:X\rightarrow Y$ is a URL-map then $X$ is complete if
and only if $Y$ is complete.
\end{proposition}

\begin{proof}
We shall prove necessity and sufficiency simultaneously. Suppose that $X$
(resp. $Y$) is complete. Let $\{y_{i}\}$ be a Cauchy sequence in $Y$ (resp. $%
X$). Choose a subsequence $\{y_{i_{j}}\}$ so that $%
d(y_{i_{j}},y_{i_{j+1}})<2^{-j}$. Let $\gamma _{j}$ be an arclength
parameterized rectifiable path in $Y$ (resp. $X$) from $y_{i_{j}}$ to $%
y_{i_{j+1}}$ such that $L(\gamma _{j})\leq 2^{-j+1}$ and $\Gamma
_{j}:=\gamma _{1}\ast \cdot \cdot \cdot \ast \gamma _{j-1}$. Take any $%
x_{1}\in f^{-1}(y_{i_{1}})$ (resp. $x_{1}=f(y_{i_{1}})$) and for $j>1$
define $x_{j}$ to be the endpoint of the lift $(\Gamma _{j})_{L}$ of $\Gamma
_{j}$ starting at $x_{1}$ (resp. of the path $f\circ \Gamma _{j}$). Since
these lifts preserve length (resp. $f$ is an arcwise isometry), it follows
by the triangle inequality that $d(x_{j},x_{j+1})\leq 2^{-j+1}$ and hence $%
\{x_{j}\}$ is Cauchy, hence convergent to a point $x\in X$ (resp. $x\in Y$).
It is clear that there exists a unique arclength parameterized path $\Gamma $
joining $x_{1}$ to $x$ such that for sufficiently large $j$, $\Gamma $%
coincides with the paths $(\Gamma _{j})_{L}$ (resp. $f\circ \Gamma _{j}$) on
any closed subsegment of the domain not containing the right endpoint. Then
the subsequence $\{y_{i_{j}}\}$, hence the sequence $y_{i}$, converges to
the endpoint of the path $f(\Gamma )$ (resp. the lift of $\Gamma $ starting
at $y_{i_{1}}$).
\end{proof}

\begin{lemma}
\label{mc} A map $f:X\rightarrow Y$ between length spaces is a URL-map if
and only if $f$ is 1-Lipshitz and for some choice of basepoints, $f$ is
basepoint preserving, each arclength parameterized rectifiable path $p$
starting at the basepoint has unique lift $p_{L}$ starting at the basepoint,
and $L(p)=L(p_{L})$.
\end{lemma}

\begin{proof}
The necessity of these conditions easily follows from the definition of
URL-map. Let us prove sufficiency. Assume first that $c$ is an arclength
parameterized rectifiable curve starting at $y\in Y$ and $x\in X$ satisfies $%
f(x)=y$. Let $k$ be a $\rho $ -curve from the basepoint in $X$ to $x$. Then $%
d:=f\circ k$ is rectifiable and since each of its initial segments has, by
assumption, a unique lift of the same length, $d$ is also arclength
parameterized and has the same length as $k$. Then $d\ast c$ is rectifiable
and arclength parameterized, so has a unique lift $(d\ast c)_{L}$ to the
basepoint in $X$. We may write $(d\ast c)_{L}=k\ast k^{\prime }$ for some
curve $k^{\prime }$. Then $k^{\prime }$ is the desired lift of $c$; $%
k^{\prime }$ must be unique since if it were not then $d\ast c$ would not
have a unique lift. Now suppose that $c:[0,a]\rightarrow Y$ is a rectifiable
path starting at $y$, and $f(x)=y$. Let $\mathcal{C}$ be the collection of
maximal (closed) intervals on which $c$ is constant. As is well-known, there
are a non-decreasing continuous function $h:[0,a]\rightarrow \lbrack 0,L(c)]$
and an arclength parameterized rectifiable path $c_{1}:[0,L(c)]\rightarrow Y$
such that $c=c_{1}\circ h$ and $h$ is strictly increasing everywhere except
on the intervals in $\mathcal{C}$. Let $d_{1}$ be the unique lift of $c_{1}$
at $x$. Define $c_{L}:[0,a]\rightarrow X$ by $c_{L}(t)=d_{1}\circ h(t)$.
Then $c_{L}$ has the same length as $d_{1}$, hence as $c_{1}$ and $c$.

Since $f$ is 1-Lipshitz, it now follows from the lifting property proved
above that if $c$ is rectifiable in $X$ then $f\circ c$ has the same length
as $c$. If $c$ is not rectifiable then $f\circ c$ cannot be rectifiable
either, for if it were, $f\circ c$ would have a rectifiable lift and a
non-rectifiable lift.
\end{proof}

\begin{lemma}
\label{ws} Every URL-map $f:X\rightarrow Y$ is a weak submetry. If $Y$ is a
geodesic space then $f$ is a submetry.
\end{lemma}

\begin{proof}
Since $Y$ is a length space, for any $\varepsilon >0$ and $x,y\in Y$ we may
join $x,y$ by a rectifiable path $c$ with the length less than $%
d(x,y)+\varepsilon $. By definition, $c$ has a lift $c_{L}$ of the same
length with endpoints $w,z$ such that $f(w)=x$ and $f(z)=y$. Since $X$ is a
length space, $d(w,z)\leq L(c_{L})=L(c)\leq d(x,y)+\varepsilon $ and the
proof of the first part is complete. Let suppose that $Y$ is a geodesic
space, $x$ is a fixed, and $y$ is arbitrary points in $Y.$ Then the points $%
x,y$ can be joined by a (shortest) geodesic $\gamma$. For any point $z\in
f^{-1}(x),$ there exists a unique lift $\overline{\gamma}$ of $\gamma,$
starting at $z,$ having the same length as $\gamma$ and some endpoint $w.$
Then $f(w)=y,$ $d(z,w)\leq L(\overline{\gamma})= L(\gamma)=d(x,y).$ On the
other hand, it follows from the first part that $f$ is 1-Lipschitz, so $%
d(x,y)\leq d(z,w)$ and $d(x,y)=d(z,w).$ This implies that $f$ is a submetry.
\end{proof}

\begin{proof}[Proof of Theorem \protect\ref{m1}]
The existence of $\overline{X}$ and the fact that $\overline{\phi }$ is a
weak submetry, hence $1$-Lipschitz, follow from Theorem \ref{ut}. To show
that $\overline{\phi }$ is URL, suppose first that $c:[0,L]\rightarrow X$ is
a $\rho $-curve at the basepoint. Let $\gamma _{c}$ denote the unique
geodesic in $\overline{X}$ from the basepoint $\ast $ to $c$. The following
formula is easy to check:%
\begin{equation}
c=\overline{\phi }\circ \gamma _{c}  \label{niceone}
\end{equation}%
In particular, $\gamma _{c}$ is a lift of $c$ to $\overline{X}$ starting at
the basepoint, and, being a geodesic, is a $\rho $-path. Since both curves
are arclength parameterized, $L(c)=L(\overline{\phi }\circ \gamma _{c})$.
Because $c$ is weakly normal, any lift $c_{L}$ of $c$ also has to be weakly
normal and hence, by Corollary \ref{runiq}, the unique geodesic $\gamma _{d}$
joining $\ast $ to some $d\in \overline{X}$. But Formula (\ref{niceone})
implies that $c=d$ and hence $c_{L}=\gamma _{c}$.

Now suppose that $d$ is a $\rho $-curve starting at $z\in X$ and $k\in 
\overline{X}$ satisfies $\overline{\phi }(k)=z$; that is, $k=\overline{\phi }%
\circ \gamma _{k}$ is a $\rho $-curve from the basepoint to $z$. By what we
just proved, $c:=k\star d$ has a unique lift $c_{L}$ to $\overline{X}$
starting at the basepoint. Now it is possible that there is some maximal
final segment $b$ of $k$ that, with reversed orientation, coincides with
some initial segment of $d$; let $c^{\prime \prime }$ be the remaining
segment of $d$. Finally, let $s$ be the final segment of $\gamma _{k}$ such
that $\overline{\phi }\circ s=b$ and $h$ be the final segment of $c_{L}$
that maps onto $c^{\prime \prime }$. Then it is easy to check that $%
d_{L}:=s^{-1}\ast h$ is the unique lift of $d$ starting at $k$, and $%
L(d_{L})=L(d)$.

By Lemma \ref{mc}, to finish the proof that $\overline{\phi }$ is URL we
need only consider an arbitrary arclength parameterized rectifiable curve $%
c:[0,L]\rightarrow X$ starting at the basepoint. By Corollary \ref{rep} we
have, for every $0\leq s<t\leq a$, a unique $\rho $-path $\rho
_{s,t}:[s,t]\rightarrow C:=c([0,L])$ such that $\rho _{s,t}$ is
fixed-endpoint homotopic to $c_{s,t}:=c\mid _{\lbrack s,t]}$ and $L(\rho
_{s,t})\leq L(c_{s,t})=t-s$. Define $c_{L}(t)$ to be the endpoint of the
unique lift $\rho _{t}^{L}$ of $\rho _{0,t}$ to $\overline{X}$ starting at
the basepoint, which exists by the special case proved above. Obviously $%
c_{L}$ is a lift of $c$. By uniqueness for any $s<t$ we must have $\rho
_{0,t}=\rho _{0,s}\star \rho _{s,t}$. Since we have shown above that $\rho
_{s,t}$ has a unique lift to $c_{L}(s)$ of the same length as $\rho _{s,t}$,
which by uniqueness must end at $c_{L}(t)$, we have 
\begin{equation}
d(c_{L}(s),c_{L}(t))\leq L(\rho _{s,t})\leq t-s\text{.}
\end{equation}

In other words, $c_{L}:[0,L]\rightarrow \overline{X}$ is 1-Lipschitz and $%
L(c_{L})\leq L$. Since we already know that $\overline{\phi }$ is distance
non-increasing, it follows that $L(c_{L})\geq L(c)=L$. Therefore $c_{L}$ is
a length-preserving lift of $c$. Now suppose that $d$ is any lift of $c$
starting at the basepoint. For any $t$ let $\gamma _{t}$ be the geodesic
from the basepoint to $d(t)$, which according to Corollary \ref{runiq} is
contained in the image of $d$ and is fixed-endpoint homotopic to $d$. Now $%
\overline{\phi }\circ \gamma _{t}$ is the $\rho $-curve $d(t)\in \overline{X}
$, which lies in the image of $\overline{\phi }\circ d=c$ and is
fixed-endpoint homotopic to $c$. By uniqueness, $\overline{\phi }\circ
\gamma _{t}=\rho _{0,t}$. Therefore $c_{L}(t)=d(t)$ for all $t$ and the
proof that $\overline{\phi }$ is URL is finished.

Now Part (1) of the theorem follows from Proposition \ref{complete}.

The uniqueness of $\overline{X}$ will follow from Lemma \ref{ws} once we
have proved the second part of the theorem. Given a URL-map $f:Z\rightarrow
X $, with some choice of basepoints, define $\overline{f}(c)$ to be the
endpoint of the unique lift of $c$ starting at the basepoint in $Z$.
Obviously $f\circ \overline{f}=\overline{\phi }$. For $c,k\in \overline{X}$,
the lift of $c\star k^{-1}$ is a path joining $f(c)$ and $f(k)$ having the
same length as $c\star k^{-1}=d(c,k)$, and therefore $\overline{f}$ is $1$%
-Lipschitz. Now let $\gamma $ be a rectifiable path starting at the
basepoint in $Z$. Then $f\circ \gamma $ is rectifiable in $X$, so has a lift 
$(f\circ \gamma )_{L}$ at the basepoint in $\overline{X}$. Now $\overline{f}%
\circ (f\circ \gamma )_{L}$ is a lift of $f\circ \gamma $ starting at the
basepint in $Z$ and so must be equal to $\gamma $. That is, $(f\circ \gamma
)_{L}$ is a lift of $\gamma $ starting at the basepoint in $\overline{X}$
having the same length as $\gamma $. Suppose $k$ is any lift of $\gamma $
starting at the basepoint in $\overline{X}$. Then $\gamma $ is also a lift
of $f\circ \gamma $ to $\overline{X}$ and therefore $k=(f\circ \gamma )_{L}$%
. We have checked the conditions of Lemma \ref{mc} to show that $\overline{f}
$ is a URL-map. Finally, suppose we have a URL-map $h$ that preserves the
basepoints with $f\circ h=\overline{\phi }$. For any $c\in \overline{X}$, $%
h\circ \gamma _{c}$ is a rectifiable path from the basepoint to $h(c)$,
which is also a lift of $\overline{\phi }\circ \gamma _{c}=c$ starting at
the basepoint in $Z$. Since this lift is unique, $h(c)=\overline{f}(c)$.
\end{proof}

\section{Continua, Fractals, and Manifolds}

\begin{proposition}
\label{lite}For any length space, $\overline{\phi }$ is light.
\end{proposition}

\begin{proof}
As is well-known, any connected subset of an $\mathbb{R}$-tree is arcwise
connected. Now if $\overline{\phi }^{-1}(p)$ contained a connected subset
with more than one point then that subset would contain a geodesic with
constant image, a contradiction to uniqueness of lifting of rectifiable
paths.
\end{proof}

In \cite{MNO} an $\mathbb{R}$-tree $T$ was defined to be $\mu $-universal
for a cardinal number $\mu $ if every $\mathbb{R}$-tree of valency at most $%
\mu $ isometrically embeds in $T$. Recall that for a point $t\in T$ the
valency at $T$ is the cardinality of the set of connected components of $%
T\backslash \{t\}$, and $T$ is said to have valency at most $\mu $ if the
valency of every point in $T$ is at most $\mu $. The existence, uniqueness
(up to isometry), and homogeneity of universal $\mathbb{R}$-trees for a
given cardinal number $\mu $ was proved in \cite{MNO}. In fact, the $\mu $%
-universal $\mathbb{R}$-tree is uniquely the complete $\mathbb{R}$-tree with
valency $\mu $ at each point. In \cite{DP} a more explicit construction was
given for the universal $\mu $-universal $\mathbb{R}$-tree $A_{\mu }$ for
arbitrary cardinal number $\mu \geq 2$.

\begin{definition}
Define $\rho $-paths $\alpha :[0,a]\rightarrow X$ and $\beta
:[0,b]\rightarrow X$ starting at $p$ to be equivalent if $\alpha ,\beta $
coincide on $[0,\varepsilon )$ for some $\varepsilon >0$. We denote the
cardinality of the set of the resulting equivalence classes by $\kappa _{p}$.
\end{definition}

\begin{proposition}
\label{val} Let $X$ be a length space, and $p\in X$. The valency of $%
\overline{X}$ at any point $\overline{p}\in \overline{\phi }^{-1}(p)$ is
equal to $\kappa _{p}$. If $X$ is separable then $\kappa _{p}\leq \mathfrak{c%
}=2^{\aleph _{0}}$.
\end{proposition}

\begin{proof}
The construction of $\overline{X}$ immediately implies the first statement.
If $X$ is separable, then $X$ itself has cardinality $\mathfrak{c}$ (unless
it is a point). Since every path is determined by its value at rational
numbers in its domain, then the cardinality of $\kappa _{p}$ is at most $%
(2^{\aleph _{0}})^{\aleph _{0}}=2^{\aleph _{0}\cdot \aleph _{0}}=2^{\aleph
_{0}}=\mathfrak{c}$.
\end{proof}

\begin{example}
If we consider the space $B$ consisting of countably many circles of radius $%
1$ attached at a point $p$ then there is a natural geodesic metric that is
uniformly one-dimensional, even though the space is not compact. The valency
of points in $\overline{B}$ is either $\aleph _{0}$ or $2$ depending on
whether they are in $\overline{\phi }^{-1}(p)$ or not.
\end{example}

The algebraic topology of one-dimensional continua, which are known to be $%
K(\pi ,1)$ spaces (cf. \cite{CF2}), has been studied since the 1950's. See,
for example, \cite{CC} for a good bibliography, recent work and open
questions concerning the fundamental groups of such spaces. A basic example
to consider is the Hawaiian earring $H$, metrized so that it consists of
countably many circles $\left\{ C_{i}\right\} _{i\in \mathbb{N}}$ each of
length $2^{1-i}$, all attached at a common basepoint $\ast $, and given the
induced geodesic metric.

\begin{proposition}
\label{haw} In the Hawaiian earring $H$ we have $\kappa _{p}=\mathfrak{c}$
for the point $p=\ast $ and $\kappa _{p}=2$ for any point $p\neq \ast $.
\end{proposition}

\begin{proof}
The last statement is evident. It is easy to find $\mathfrak{c}$ unit weakly
normal loops starting (and ending) at $p=\ast $ such that no two coincide on
any interval $[0,\varepsilon )$. In fact, one can define a path that wraps
one of two ways around $C_{1}$, then one of two ways around $C_{2}$, and so
on. Then reverse the parametrization, so that $C_{1}$ is wrapped around
last. It is clear that any such path is encoded as a sequence $%
\{x_{n}\}_{n\in \mathbb{N}}$ with values in the set $\{-1,1\}$. An arclength
parameterized path will be equivalent to one of these if and only if it
wraps the same way around $C_{i}$ for all sufficiently large $i$, or in
other words, defines an equivalent sequence $\{y_{n}\}$--we mean here that $%
x $ is equivalent to $y$ if there is a number $m\in \mathbb{N}$ such that $%
x_{n}=y_{n}$ for all $n\geq m$. The following lemma finishes the proof.
\end{proof}

\begin{lemma}
In the Hawaiian earring $H$, we have $\mathfrak{c}$ different equivalence
classes of the above type of sequence.
\end{lemma}

\begin{proof}
Let us take first the sequence $x_{n}\equiv 1$. Consider now an arbitrary
sequence $z=\{z_{n}\}$ of natural numbers and define subsequently two
sequences $s(z)_{n}:=\sum_{i=1}^{n}z_{i}$ and $\sigma
(z)_{n}:=\sum_{i=1}^{n}2^{s(z)_{i}}$. It is clear that if for another such
sequence $w=\{w_{n},\}$, $w_{m}\neq z_{m}$, then $\sigma (w)_{n}\neq \sigma
(z)_{n}$ for all $n\geq m$. Then for every sequence $z$ above, define
another sequence $x(z)_{n}$ with values in $\{-1,1\}$ by the equations $%
x(z)_{m}=-1$ if $m=\sigma (z)_{k}$ for some $k\in \mathbb{N}$, and $%
x(z)_{m}=1$ for all other $m\in \mathbb{N}$. Then it follows from the
statement above that $x(z)_{n}$ is not equivalent to $x$ and is not
equivalent to $x(w)$ if $w\neq z$. So we get $\aleph _{0}^{\aleph _{0}}=%
\mathfrak{c}$ pairwise non-equivalent sequences of above sort.
\end{proof}

\begin{proof}[Proof of Theorem \protect\ref{universal}]
The first part of the theorem is an immediate consequence of Proposition \ref%
{val} and the theorem from \cite{MNO} that every $\mathbb{R}$-tree of
valency at most $\mathfrak{c}$ isometrically embeds in $A_{\mathfrak{c}}$.
Next, for any point $p\in S_{c}$ there is clearly a bi-Lipschitz embedding $%
h:H\rightarrow S_{c}$ such that $h(\ast )=p$. Therefore $\kappa _{p}\geq 
\mathfrak{c}$. The case of $S_{g}$ is more tricky. There are countably many
rectifiable loops $C_{i}$ in $S_{g}$ starting at any fixed point $p$ such
that every $C_{i}$ is a topologically embedded circle and $L(C_{i+1})<\frac{1%
}{3}L(C_{i})$ for all natural numbers $i$. Then one can prove with a little
more details than in Proposition \ref{haw}, that $\kappa _{p}\geq \mathfrak{c%
}$. A similar argument now finishes the proof of the theorem.
\end{proof}

\section{URL($G$)-maps}

Throughout this section, let $X$ and $Y$ be length spaces. The following
statement easily follows from definitions.

\begin{proposition}
\label{cat} The collection of pointed length spaces with URLs as morphisms
is a category, which we will refer to as the URL category.
\end{proposition}

\begin{proposition}
\label{samet} Let $f:(X,\ast )\rightarrow (Y,\ast )$ be a basepoint
preserving URL-map of length spaces. Then there is a commutative diagram 
\begin{equation*}
\begin{array}{lll}
(\overline{X},\ast ) & \overset{\overline{f}}{\longrightarrow } & (\overline{%
Y},\ast ) \\ 
\downarrow ^{\overline{\phi _{X}}} &  & \downarrow ^{\overline{\phi _{Y}}}
\\ 
(X,\ast ) & \overset{f}{\longrightarrow } & (Y,\ast ) \\ 
&  & 
\end{array}%
\end{equation*}%
of URL-maps preserving basepoints, with unique $\overline{f}$, where $%
\overline{\phi _{X}}$ and $\overline{\phi _{Y}}$ are $\mathbb{R}$-tree
covering maps for $X$ and $Y$ respectively, and $\overline{f}$ is isometry.
\end{proposition}

\begin{proof}
By Theorem \ref{m1}, there are unique basepoint preserving URL-maps $\psi :(%
\overline{Y},\ast )\rightarrow (X,\ast )$ such that $f\circ \psi =\overline{%
\phi _{Y}},$ and $\overline{f}:(\overline{X},\ast )\rightarrow (\overline{Y}%
,\ast )$ such that $\psi \circ \overline{f}=\overline{\phi _{X}}$. Then $%
\overline{\phi _{Y}}\circ \overline{f}=f\circ \psi \circ \overline{f}=f\circ 
\overline{\phi _{X}}$. By Proposition \ref{cat}, the composition $f\circ 
\overline{\phi _{X}}$ is a URL-map. It follows from Theorem \ref{m1} that
there exists a unique URL-map $\overline{g}:(\overline{Y},\ast )\rightarrow (%
\overline{X},\ast )$ such that $f\circ \overline{\phi _{X}}\circ \overline{g}%
=\overline{\phi _{Y}}$. By the previous argument, $\overline{\phi _{X}}\circ 
\overline{g}=\psi $, and $\overline{\phi _{X}}\circ \overline{g}\circ 
\overline{f}=\psi \circ \overline{f}=\overline{\phi _{X}}$. Now Theorem 1
and the previous equations imply that $\overline{g}\circ \overline{f}=1_{X}$%
. By a similar argument, $\overline{\phi _{Y}}\circ (\overline{f}\circ 
\overline{g})=f\circ \overline{\phi _{X}}\circ \overline{g}=f\circ \psi =%
\overline{\phi _{Y}}$, and $\overline{f}\circ \overline{g}=1_{Y}$. Thus $%
\overline{f}$ and $\overline{g}$ are inverses ofone another, and both weak
submetries (hence, 1-Lipschitz maps) by Proposition \ref{ws}. Then both $%
\overline{f}$ and $\overline{g}$ are (surjective) isometries.
\end{proof}

So in the situation of Proposition \ref{samet}, $(\overline{X},\ast )$ can
be identified with $(\overline{Y},\ast )$.

\begin{lemma}
\label{rhl}If $f:X\rightarrow Y$ is URL then a path $c$ in $X$ is a $\rho $%
-path if and only if $f\circ c$ is a $\rho $-path.
\end{lemma}

\begin{proof}
Identifying $(\overline{Y},\ast )$ with $(\overline{X},\ast )$ via the
isometry $\overline{f}^{-1}$, we get that $f\circ \overline{\phi _{X}}=%
\overline{\phi _{Y}}$ (see Proposition \ref{samet}). This implies that $c$
and $f\circ c$ have the same lift $\gamma _{c}$ to $\overline{X}$ with
respect to the URL-maps $\overline{\phi _{X}}$ and $\overline{\phi _{Y}}$.
As in the proof of Theorem \ref{m1} either of the curves $c$ or $f\circ c$
is $\rho $-curve if and only if $\gamma _{c}$ is a geodesic in $\overline{X}$%
. This finishes the proof.
\end{proof}

\begin{proposition}
\label{inje} If $f:X\rightarrow Y$ is URL then letting $f_{\ast }(c):=f\circ
c$ defines an injective homomorphism $f_{\ast }:\Lambda (X)\rightarrow
\Lambda (Y)$. In addition, the orbits of the subgroup $f_{\ast }(\Lambda
(X)) $ in $\overline{Y}$ are closed, and $\phi _{X}$ and $f$ are naturally
identified, respectively, with the metric quotient map $\psi :\overline{Y}%
\rightarrow f_{\ast }(\Lambda (X))\backslash \overline{Y}$ and the map $\phi
:f_{\ast }(\Lambda (X))\backslash \overline{Y}\rightarrow \Lambda
(Y)\backslash \overline{Y}$ induced by the inclusion map $f_{\ast }(\Lambda
(X))\rightarrow \Lambda (Y)$.
\end{proposition}

\begin{proof}
If $c$ is the cancelled concatenation of $c_{1}$ and $c_{2}$ in $\Lambda (X)$%
, then by uniqueness $f\circ c$ is the cancelled concatenation of $f\circ
c_{1}$ and $f\circ c_{2}$, which is precisely what it means for $f_{\ast
}:\Lambda (X)\rightarrow \Lambda (Y)$ to be a homomorphism. The injectivity
of $f_{\ast }$ immediately follows from definitions and Lemma \ref{rhl}. Now
if we identify $\overline{X}$ with $\overline{Y}$ via the isometry $%
\overline{f}$, then $\Lambda (X)$ is naturally identified with $f_{\ast
}(\Lambda (X))$ while the lift of any element $c\in \Lambda (X)$ with
respect to $\overline{\phi _{X}}:\overline{X}\rightarrow X$ is identified
with the lift of $f\circ c$ with respect to $\overline{\phi _{Y}}$, because $%
f\circ \overline{\phi _{X}}=\overline{\phi _{Y}}$. This, together with Lemma %
\ref{ws}, implies the other statements of the proposition.
\end{proof}

\begin{definition}
\label{URL}A URL-map $f:(X,\ast )\rightarrow (Y,\ast )$ is called URL($G$)
for some subgroup $G$ of $\Lambda (Y)$ if for every $c\in G$, its lift $%
c_{L} $ to $X$ starting at $\ast $ is a loop.
\end{definition}

If $H$ is a subgroup of $G$ then obviously every URL($G$)-map is a URL($H$%
)-map. By Lemma \ref{rhl}, the lift $c_{L}$ in the above definition is in
fact a $\rho $-loop. Proposition \ref{inje} immediately implies:

\begin{corollary}
\label{gloop}If $f:X\rightarrow Y$ is URL then:

\begin{enumerate}
\item If $G$ is a subgroup of $\Lambda (Y)$ and $f$ is URL($G$) then $%
f^{\ast }(G):=\{c_{L}:c\in G\}$ is a subgroup of $\Lambda (X)$, isomorphic
to $G$.

\item If $H$ is a subgroup of $\Lambda (X)$ then $f_{\ast }(H):=\{f\circ
c\in \Lambda (Y):c\in H\}$ is a subgroup of $\Lambda (Y)$, isomorphic to $H$.
\end{enumerate}
\end{corollary}

\begin{example}
Suppose $f:(X,\ast )\rightarrow (Y,\ast )$ is a traditional covering map, $X$
is connected, $Y$ is a length space, and $X$ has the lifted metric. Since
traditional metric covering maps are URL-maps and fibrations, $\pi $ is URL($%
H_{T}$), where $H_{T}$ is the group of null-homotopic $\rho $-loops in $Y$
at $\ast $.
\end{example}

\begin{definition}
\label{gsimp}If $G$ is a subgroup of $\Lambda (Y)$ then $X$ is called $G$%
-universal if there is a URL($G$) function $f:X\rightarrow Y$ such that if $c
$ is a $\rho $-loop at $\ast $ in $X$ then $f\circ c\in G$. In this case the
map $f$ will also be called $G$-universal.
\end{definition}

From Proposition \ref{inje} we have:

\begin{proposition}
\label{gloopy}A URL-map $f:X\rightarrow Y$ is URL($G$) if and only if $%
G\subset f_{\ast }(\Lambda (X))$, and $G$-universal if and only if $f_{\ast
}(\Lambda (X))=G$.
\end{proposition}

\begin{proposition}
\label{fund}The map $h_{\Lambda }:\Lambda (X)\rightarrow \pi _{1}(X)$ that
takes each $\rho $-loop to its homotopy equivalence class, is a homomorphism
with the kernel $H_{T}$ and image $\mu _{1}(X)$, the subgroup of $\pi _{1}(X)
$ consisting of all equivalence classes having at least one rectifiable
representative. Therefore, if each element of the fundamental group of $X$
contains a rectifiable representative, then $\pi _{1}(X)$ is naturally
isomorphic to $\Lambda (X)/H_{T}:=\pi _{1}^{H_{T}}(X)$.
\end{proposition}

\begin{proof}
The map $h_{\Lambda }$ is homomorphism, since the canceled concatenation of
two $\rho $-curves is fixed-endpoint homotopic to their concatenation.
Evidently, its kernel is $H_{T}$; the statement about image follows from
Corollary \ref{rep}.
\end{proof}

\begin{remark}
It is not true in general, even for compact geodesic spaces $X$, that each
element of the fundamental group of $X$ contains a rectifiable
representative (cf. \cite{BPS}), contrary to Remark 1.13 (b) on p. 10 in 
\cite{Gr}.
\end{remark}

\begin{proof}[Proof of Corollary \protect\ref{slsc}]
Let $\pi :Y\rightarrow X$ be the traditional universal covering map.
Moreover, since $Y$ is simply connected, the image of every loop in $Y$ is
an element of $H_{T}$, and this means $\pi $ is $H_{T}$-universal. From
Proposition \ref{inje} it follows that $Y$ is naturally isometric to ${H_{T}}%
\backslash \overline{X}$. It is well-known (and easy to check) that since $X$
is semi-locally simply connected, every path contains a piecewise geodesic
in its fixed-endpoint homotopy class, so the second part of the corollary
follows from Proposition \ref{fund}.
\end{proof}

We will need the following lemma, which is a kind of metric Second
Isomorphism Theorem.

\begin{lemma}
\label{metbasic}Let $G$ act by isometries on $X$ and $f:X\rightarrow
Y=G\backslash X$ be the quotient map. Suppose $H$ is a subgroup of $G$ and $%
\psi :X\rightarrow H\backslash X=Z$ is the quotient map.

\begin{enumerate}
\item There is a unique function $\phi :Z\rightarrow Y$ such that the
following diagram commutes: 
\begin{equation*}
\begin{array}{lll}
X & \overset{\psi }{\longrightarrow } & \text{ \ \ \ \ \ \ }Z \\ 
\text{ \ \ \ }\searrow ^{f} &  & ^{\phi }\swarrow \\ 
& Y & 
\end{array}%
\end{equation*}

\item If $H$ has closed orbits and $Z$ is given the quotient metric then $%
\phi $ is a weak submetry.

\item If $H$ is a normal subgroup of $G$ with closed orbits then $G/H$ acts
naturally as isometries on $Z$ and $\phi $ is the quotient mapping onto $%
Y=\left( G/H\right) \backslash Z$.
\end{enumerate}
\end{lemma}

\begin{proof}
Define $\phi (\psi (x))=f(x)$; this is the only possible way to define $\phi 
$ to make the diagram commute, so we need only check that it is
well-defined. But if $\psi (x)=\psi (y)$ then by definition both $x$ and $y$
lie in the same orbit of $H$, hence the same orbit of $G$, so $f(x)=f(y)$.
For the second part note that $\phi $ is distance non-increasing because the
orbits of $H$ are contained in the orbits of $G$. Let $x,y\in Y$ and $%
\varepsilon >0$. Since $f$ is a weak submetry there are $a\in f^{-1}(x),b\in
f^{-1}(y)$ such that $d(a,b)<d(x,y)+\varepsilon $. Then by definition $%
d(Ha,Hb)\leq d(a,b)<d(x,y)+\varepsilon $, $\phi (Ha)=x$ and $\phi (Hb)=y$,
which proves that $\phi $ is a weak submetry.

For the third part define $gH(Hx)=Hg(x)$ for $x\in X$ and $g\in G$. Since $H$
is normal it is easy to check that we have a well-defined function
corresponding to the coset $gH\in G/H$. Given $g_{1},g_{2}\in G$, 
\begin{equation*}
\left( g_{1}g_{2}\right) H(Hx)=Hg_{1}(g_{2}(x))=g_{1}H(Hg_{2}(x))
\end{equation*}%
\begin{equation*}
=g_{1}H(g_{2}H(Hx))=(g_{1}H\circ g_{2}H)(Hx)\text{.}
\end{equation*}%
That is, we have a properly defined action. Note that this also implies that
each function $gH$ has an inverse and hence is injective.

Moreover, since $d(Hx,Hy)=\inf \{d(k(x),y):k\in H\}$ and $g$ is an isometry, 
\begin{equation*}
d(Hg(x),Hg(y))=\inf \{d(kg(x),g(y)):k\in H\}
\end{equation*}%
\begin{equation*}
=\inf \{d(g^{-1}kg(x),y):k\in H\}=d(Hx,Hy)
\end{equation*}%
Now $Hx,Hy$ lie in the same orbit of this action if and only if $Hg(x)=Hy$
for some $g\in G$; equivalently $x$ and $y$ lie in the same orbit of $G$,
which is equivalent to $f(x)=f(y)$. By definition this is equivalent to $%
\phi (Hx)=\phi (Hy)$. This shows that the orbits of the action of $G/H$ are
the same as the inverse images of $\phi $. Finally we need to check that $Y$
has the quotient metric with respect to this action. But the distance
between the orbits of $Hx$ and $Hy$ is the infimum of the distances between $%
Hx$ and $Hg(y)$, where $g\in G$. This in turn is the infinum of distances
between $x$ and $kg(y)$, where $k\in H$, which is the distance between the
orbits of $x$ and $y$ with respect to $G$. But since $Y$ has the quotient
metric, this is precisely $d(f(x),f(y))=d(\phi (Hx),\phi (Hy))$.
\end{proof}

\begin{notation}
\label{added}For any subgroup $G$ of $\Lambda (Y)$ that has closed orbits in 
$\overline{Y}$ we will denote $G\backslash \overline{Y}$ with the quotient
metric by $\overline{Y}^{G}$. We will use the notation $\psi _{G}:\overline{Y%
}\rightarrow G\backslash \overline{Y}$ for the quotient mapping and $%
\overline{\phi }^{G}:\overline{Y}^{G}\rightarrow Y=\Lambda (Y)\backslash 
\overline{Y}$ for the mapping analogous to $\phi $ from Lemma \ref{metbasic}.
\end{notation}

\begin{proposition}
\label{induced}Let $G$ be a subgroup of $\Lambda (Y)$ and suppose $%
f:X\rightarrow Y$ is a URL($G$)-map. Then

\begin{enumerate}
\item There is a unique (up to basepoint choice) function $g:\overline{Y}%
^{G}\rightarrow X$ such that $\overline{\phi }^{G}=f\circ g$, and moreover $%
g $ is a weak submetry if $G$ has closed orbits in $\overline{Y}$.

\item If $f$ is $G$-universal then $G$ has closed orbits in $\overline{Y}$, $%
g$ is an isometry, and $f$ can be identified with $\phi_G$.
\end{enumerate}
\end{proposition}

\begin{proof}
Using Proposition \ref{samet} and Proposition \ref{inje} we identify $%
\Lambda (X)$ with $f_{\ast }(\Lambda (X))\subset \Lambda (Y)$ and $\phi _{X}$
with $\psi :\overline{Y}\rightarrow f_{\ast }(\Lambda (X))\backslash 
\overline{Y}$. Proposition \ref{gloopy} implies that $G\subset f_{\ast
}(\Lambda (X))$, and from Lemma \ref{metbasic} we have a unique (up to
basepoint choice) mapping $g:\overline{Y}^{G}\rightarrow f_{\ast }(\Lambda
(X))\backslash \overline{Y}=X$, which is a weak submetry when the orbits of $%
G$ are closed. With these identifications, $f\circ g=\phi _{G}$, and
uniqueness follows from the fact that these identifications are all uniquely
determined by the basepoint choice.

The second part follows from Proposition \ref{gloopy}, since in this case $%
G=f_{\ast }(\Lambda (X))$, $g:\overline{Y}^{G}\rightarrow f_{\ast }(\Lambda
(X))\backslash \overline{Y}$ is the identity map, and under our
identifications $f$ is identified with $\phi _{G}$.
\end{proof}

URL($G$)-maps have the following useful property:

\begin{proposition}
\label{diag2}Let $X,Y,Z$ be length spaces and assume that the following
diagram of continuous base-preserving maps is commutative: 
\begin{equation*}
\begin{array}{lll}
X & \overset{h}{\longrightarrow } & \text{ \ \ \ \ \ \ }Z \\ 
\text{ \ \ \ }\searrow ^{f} &  & ^{g}\swarrow  \\ 
& Y & 
\end{array}%
\end{equation*}%
If $g$ is a URL-map, then $f$ is a URL-map if and only if $h$ is a URL-map,
and each of the maps $f$ or $h$ is uniquely defined by the another one and $g
$.

Suppose further that $G$ is a subgroup of $\Lambda (Y,\ast )$, and $g$ is
URL($G$).

\begin{enumerate}
\item If $h$ is URL($g^{\ast }(G)$) then $f$ is URL($G$).

\item If $f$ is URL($G$) then $h$ is a uniquely determined (up to basepoint
choice) URL($g^{\ast }(G)$)-map.
\end{enumerate}
\end{proposition}

\begin{proof}
In checking that maps are URL we will use Lemma \ref{mc} without further
reference, and \textit{all paths will be assumed to start at the basepoint}.
Suppose first that $h$ is URL. Then $f=g\circ h$ is URL by Proposition \ref%
{cat}. We must have $f=g\circ h,$ so $f$ is uniquely determined by $h$ and $%
g $.

Now suppose that $f$ is URL and $c$ is a path in $Z$ of the length $L<\infty 
$. Then $c_{Y}:=g\circ c$ has length $L$, and thus has unique lift $c_{X}$
to $X$ of length $L$. Then $g\circ h\circ c_{X}=f\circ c_{X}=c_{Y}$ and $%
g\circ c=c_{Y}$, so $c=h\circ c_{X}$ because $g$ is URL. Thus $c_{X}$ is a
lift of $c$. If there is another such lift $c_{X}^{\prime }$, then $f\circ
c_{X}^{\prime }=g\circ h\circ c_{X}^{\prime }=g\circ c=c_{Y}$ and $%
c_{X}^{\prime }=c_{X}$, because $c_{Y}$ is a rectifiable path and $f$ is
URL. Moreover, $L(h\circ c_{X})=L(g\circ (h\circ c_{X}))=L(f\circ
c_{X})=L(c) $ because $g$ and $f$ are URL. Thus $h$ preserves the length of
paths, and $h $ is URL. Let $x\in X$ and $\gamma $ be any $\rho $-curve from 
$\ast $ to $x$ in $X$. Then $h\circ \gamma $ must be the unique lift of $%
f\circ \gamma $ to $Z$ and therefore $h(x)$ is uniquely determined by $f$
and $g$ (and the basepoints).

Suppose now that $G$ is a subgroup of $\Lambda (Y,\ast )$, and $g$ is URL($G$%
).

(1) If $h$ is URL($g^{\ast }(G)$) and $c\in G$ then since $g$ is URL($G$), $%
c_{L}\in g^{\ast }(G)$ ($c_{L}$ is the lift of $c$ to $Z$), and hence lifts
as a loop $d$ in $X$. However, by a similar argument to the above, $d$ is
the unique lift of $c$ to $X$ (we have already shown $f$ is URL!) and this
shows that $f$ is URL($G$).

(2) If $f$ is URL($G$) and $c\in g^{\ast }(G)$ then the lift $d$ of $g\circ
c\in G$ to $X$ is a loop. Now $g\circ h\circ d=f\circ d=g\circ c$ and hence $%
h\circ d$ is the unique lift of $g\circ c$ to $Z$. But then $c=h\circ d$.
That is, the loop $d$ is the unique lift of $c$ to $X$, which proves that $h$
is URL($g^{\ast }(G)$).
\end{proof}

\begin{remark}
Similar arguments show that in the setting of Proposition \ref{diag2}, $g$
is URL provided both $f$ and $h$ are URL and $g$ has the property that for
any non-rectifiable path $c$ in $Z$, $g\circ c$ is non-rectifiable.
\end{remark}

\begin{proof}[Proof of Theorem \protect\ref{m2}]
The first statement follows from Lemma \ref{metbasic} and Theorem \ref{ut}.
Since $\overline{\phi }^{G}$ is URL by assumption, similar to the proof of
Proposition \ref{induced}, there exists unique mapping $\psi :\overline{X}%
\rightarrow \overline{X}^{G}$ such that $\overline{\phi }^{G}\circ \psi =%
\overline{\phi }$. Then $\psi =\psi _{G}$ and it follows from Proposition %
\ref{samet} that $\overline{X}$ and $\psi _{G}$ are naturally identified
respectively with $\overline{\overline{X}^{G}}$ and the $\mathbb{R}$-tree
universal covering map for $\overline{X}^{G}$. It is clear now by
Proposition \ref{inje} that $\overline{\phi }_{\ast }^{G}(\Lambda (\overline{%
X}^{G}))=G$, and $\overline{\phi }^{G}$ is $G$-universal by Proposition \ref%
{gloopy}.

The second part of the theorem follows by taking $f^{G}:=g$ from Proposition %
\ref{induced} (with $X$ and $Y$ interchanged) and observing that $f^{G}$ is
uniquely determined and URL by Proposition \ref{diag2}. If $f$ is $G$%
-universal, then by Proposition \ref{induced} (2), $f^{G}$ is an isometry.
This implies the uniqueness of $\overline{\phi }^{G}$ as a $G$-universal map
onto $X$, and together with the above considerations, the first part of the
theorem.

The third part is a consequence of Lemma \ref{metbasic}. The fourth part
follows from Proposition \ref{complete} because $\overline{\phi}^G$ is URL.
\end{proof}

\begin{remark}
Following the classical usage, it now makes sense to say that a URL-map $%
f:X\rightarrow Y$ is called normal (or regular) if $f_{\ast }(\Lambda (X))$
is a normal subgroup of $\Lambda (Y)$. We may then refer to $\Lambda
(Y)/f_{\ast }(\Lambda (X))$ as the group of deck transformations of $f$.
\end{remark}

\begin{example}
Let $X$ be the Euclidean plane, and let $c$ and $d$ be different semicircles
on the same circle of circumference $2$ parameterized as $\rho $-paths with
the same starting point $\ast $ and endpoint $x$. Take sequences $%
c(t_{i})\rightarrow x$ and $d(t_{i})\rightarrow x$ with $t_{i}$ strictly
increasing. Define $c_{i}:=\left( d\mid _{\lbrack 0,t_{i}]}\right) \star
f\star \left( c\mid _{\lbrack 0,t_{i}]}\right) ^{-1}$ where $f$ is an
arclength parameterized straight line from $d(t_{i})$ to $c(t_{i})$. Let $G$
be the subgroup of $\Lambda (X)$ generated by the loops $c_{i}$. Note that
no loop in $G$ passes through the point $x$, and therefore for any $k\in G$, 
$d(k,d\star c^{-1})>1$. In other words, $d\star c^{-1}$ is not in the
closure of $G$. This shows that simply being a closed subgroup of $\Lambda
(X)$ is not sufficient to have closed orbits.
\end{example}

\begin{lemma}
\label{sc}If $G$ is a subgroup of $\Lambda (X)$ then the orbits of $G$ are
closed in $\overline{X}$ if and only if $G$ is \textit{strongly closed }the
following sense: Given $c_{i}\in G$ and $c\in \overline{X}$ such that $%
c_{i}\star c\rightarrow d$ then $d\star c^{-1}\in G$.
\end{lemma}

\begin{proof}
That the orbits of $G$ are closed means precisely that if $c\in \overline{X}$
and $c_{i}\in G$ such that $c_{i}\star c\rightarrow d$, then there is some $%
a\in G$ such that $d=a\star c$. If $d\star c^{-1}\in G$ then certainly it
plays the role of $a$. On the other hand, if such an $a$ exists then since $%
a $ is a loop, $c$ and $d$ have the same endpoint. Then we may concatenate
on each side to obtain $d\star c^{-1}=a$.
\end{proof}

\begin{remark}
If $G$ is strongly closed then given a convergent sequence $\left(
c_{i}\right) $ in $G$, we may apply the strongly closed criterion with $c$
the trivial loop to see that the limit of $(c_{i})$ is actually in $G$,
implying that $G$ is closed.
\end{remark}

The next Proposition shows that for any strongly closed group $G$, the
mapping $\overline{\phi }^{G}$ is tantalizingly close to being URL.

\begin{proposition}
\label{almost URL}If $G$ is any strongly closed subgroup of $\Lambda (X)$
then $\overline{\phi }^{G}:\overline{X}^{G}\rightarrow X$ is a weak submetry
such that

\begin{enumerate}
\item For every $x,y\in \overline{X}^{G}$, $d(x,y)$ is the infimum of the
lengths of rectifiable paths $c$ such that $\overline{\phi }^{G}\circ c$ is
a $\rho $-path in $X$ of the same length as $c$.

\item For every rectifiable path $k$ starting at $p\in X$ and $q\in \left( 
\overline{\phi }^{G}\right) ^{-1}(p)$ there is some lift $d$ of $k$ to $%
\overline{X}^{G}$ of the same length as $k$ and starting at $q$.
\end{enumerate}
\end{proposition}

\begin{proof}
That $\overline{\phi }^{G}$ is a weak submetry follows from Lemma \ref%
{metbasic}. Since $\psi _{G}$ is a quotient map, for $x,y\in \overline{X}%
^{G} $\ and $\varepsilon >0$ we may find a geodesic $\gamma $ joining points 
$w,z$ with $\psi _{G}(w)=x$, $\psi _{G}(z)=y$, and $d(x,y)\leq L(\gamma
)<d(x,y)+\varepsilon $. Let $c:=\psi _{G}\circ \gamma $. Since $\psi _{G}$
and $\overline{\phi }^{G}$ are both distance (hence length) nonincreasing, 
\begin{equation*}
L(c)\geq L(\overline{\phi }^{G}\circ c)=L(\overline{\phi }\circ \gamma
)=L(\gamma )\geq L(c)
\end{equation*}%
Letting $\varepsilon \rightarrow 0$ finishes the proof of the first part.

For the second part, let $q^{\prime }\in \overline{X}$ be such that $\psi
_{G}(q^{\prime })=q$, let $\gamma _{k}$ be the lift of $k$ to $q^{\prime }$
(which is a geodesic) and let $d:=\psi _{G}\circ \gamma _{k}$. Since $%
\overline{\phi }^{G}\circ \psi _{G}=\overline{\phi }$, $d$ is a lift of $k$,
and by the same argument as in the first part, $d$ has the same length as $k$%
.
\end{proof}

A basic question remains: whether $\overline{\phi }^{G}:\overline{X}%
^{G}\rightarrow X$ from Proposition \ref{almost URL} is a URL-map for every
strongly closed subgroup of $\Lambda (X)$?

\section{The Metric Core}

In this section we construct a large family of URL-maps for spaces that may
not be semi-locally simply connected. In many cases the maps involve domains
that are simply connected, and in some cases even $CAT(0)$. For this section
the reader is referred to the appendix for notation and background. Let $%
H_{L}$ denote the group of $\rho $-paths that are null-homotopic via a $1$%
-Lipschitz homotopy. In this section let $X$ and $Y$ be length spaces.

\begin{definition}
For $x,y\in \widetilde{X}$, define 
\begin{equation*}
d(x,y):=\inf \{L(\phi \circ \gamma )\}
\end{equation*}%
where the infimum is taken over all paths $\gamma $ joining $x$ and $y$, and
we take the convention that the infimum of the empty set is $\infty $.
Define 
\begin{equation*}
\widehat{X}:=\left\{ x\in \widetilde{X}:d(x,\ast )<\infty \right\} \text{.}
\end{equation*}%
The set $\widehat{X}$ will be called the metric core of $\widetilde{X}$ and
we will denote the restriction of $\phi $ to $\widehat{X}$ by $\widehat{\phi 
}$, and the restrictions of the projections $\phi _{i}$ will be denoted by $%
\widehat{\phi _{i}}$; we will also denote $\widehat{\phi }$ by $\widehat{%
\phi }_{0}$.
\end{definition}

\begin{definition}
Let $\mu (X):=\lambda (\mu _{1})\subset \delta (X)$. We will call $\mu (X)$
the uniform metric deck group of $X$.
\end{definition}

\begin{theorem}
\label{core}The function $d$ defined above defines a (finite) length space
metric on $\widehat{X}$ such that

\begin{enumerate}
\item The inclusion of $\widehat{X}$ with this metric into $\widetilde{X}$
is uniformly continuous.

\item If $Y$ is a length space and $g:Y\rightarrow \widetilde{X}$ is a
function such that $g(Y)\cap \widehat{X}\neq \varnothing $ and $\phi \circ g$
is $1$-Lipschitz then $g(Y)\subset \widehat{X}$ and $g$ is $1$-Lipschitz
(hence continuous) in the metric of $\widehat{X}$.

\item Every $1$-Lipschitz homotopy between rectifiable paths in $X$ lifts to 
$\widehat{X}$ and in particular $\widehat{\phi }$ is an $H_{L}$-URL-map.

\item If $X$ is complete then $\widehat{X}$ is complete.

\item The group $\mu (X)$ is the stabilizer of $\widehat{X}$ in the uniform
fundamental group $\delta (X)$ and acts isometrically and freely on $%
\widehat{X}$ with metric quotient $\mu (X)\backslash\widehat{X}=X$.
\end{enumerate}
\end{theorem}

\begin{proof}
We will use Notation \ref{indexnot} from the appendix. As mentioned in the
appendix, for any $i$, $\phi _{i}:\widetilde{X}\rightarrow X_{i}$ is the
UU-covering of $X_{i}$ and therefore this mapping has the same lifting
properties as $\phi $. Symmetry, positive definiteness, and the triangle
inequality (with possibly infinite values and the usual conventions for
adding extended real numbers) of $d$ are clear from the definition. From
this it follows that if $x,y\in \widehat{X}$ then $d(x,y)<\infty $, hence $d$
is a (finite) metric on $\widehat{X}$. The fact that $\widehat{X}$ is a
length space will follow from the definition of the metric if we show that $%
\widehat{\phi }$ is length-preserving. Since the distance between any two
points in $X$ is the infimum of lengths of paths joining them, including
paths that may not be projections of paths in $\widetilde{X}$, $\widehat{%
\phi }$ is distance non-increasing and hence length non-increasing. If $%
\widehat{\phi }\circ \gamma :[0,1]\rightarrow X$ is rectifiable then for any
partition $t_{0},...,t_{k}$ of $[0,1]$ we have by definition of the distance
in $\widehat{X}$, 
\begin{equation*}
\sum_{i=1}^{k}d(\gamma (t_{i}),\gamma (t_{i-1}))\leq \sum_{i=1}^{k}L(%
\widehat{\phi }\circ \gamma \mid _{\lbrack t_{i-1},t_{i}]})=L(\widehat{\phi }%
\circ \gamma )<\infty
\end{equation*}%
and therefore $\gamma $ is rectifiable and $L(\gamma )\leq L(\widehat{\phi }%
\circ \gamma )$. Therefore $\widehat{\phi }$ is length-preserving. A similar
argument shows that every $\widehat{\phi }_{i}$ is similarly
length-preserving, a fact that we will need below.

To show that the inclusion is uniformly continuous, due to the compatibility
properties of $\widetilde{d}$ discussed in the appendix (the uniform
structure is the inverse limit structure), we need only show the following:
For any $\varepsilon >0$ there exists some $\delta >0$ such that if $x,y\in 
\widehat{X}$ satisfy $d(x,y)<\delta $ then $d(\phi _{i}(x),\phi
_{i}(y))<\varepsilon $ for any $i$. But in fact each $\phi _{i}$ is distance
non-increasing and we may simply take $\delta =\varepsilon $.

For Part (2), suppose that $g(y)\in \widehat{X}$. Then any $x\in Y$ may be
joined to $y$ via a rectifiable path $c$. Since $\phi \circ g$ is distance
non-increasing, $\phi \circ c$ is rectifiable, which shows that the distance
between $g(x)$ and $g(y)$ is finite, i.e. $g(x)\in \widehat{X}$. Choosing
the length of $c$ to be close to $d_{Y}(x,y)$ and applying the definition of
the distance in $\widehat{X}$ shows that $g$ is distance non-increasing into 
$\widehat{X}$.

For Part (3), let $c$ be a $\rho $-path in $X$. Then the unique lift of $c$
to the basepoint (which lies in $\widehat{X}$) satisfies the conditions of
the second part of this theorem, and so the unique lift of $c$ to $%
\widetilde{X}$ to the basepoint must be a path in $\widehat{X}$. We have
already observed that $\widehat{\phi }$ is length preserving, so the lift
has the same length. If there were another lift of $c$ to $\widehat{X}$ then
since the inclusion of $\widehat{X}$ into $\widetilde{X}$ is continuous
there would be two lifts of $c$ in $\widetilde{X}$, a contradiction, showing
that $\widehat{\phi }$ is a URL-map. Similarly, any $1$-Lipschitz homotopy
at the basepoint lifts uniquely to $\widehat{X}$, finishing the proof of the
third part.

For the fourth part, let $\left( z_{j}\right) _{j=1}^{\infty }$ be a Cauchy
sequence in $\widehat{X}$. Letting $x_{j}:=\phi (z_{j})\in X$, the fact that 
$\phi $ is distance non-increasing implies that $\{x_{j}\}$ is Cauchy with
limit $x$. Let $\left\{ \gamma _{j}\right\} $ be a collection of paths such
that $\gamma _{j}$ joins $x_{j}$ to $x$ and $L(\gamma _{j})<2d(x_{j},x)$.
Define $\alpha _{j}$ to be the unique lift of $\gamma _{j}$ starting at $%
z_{j}$. Since $\phi $ is length preserving, if $y_{j}$ denotes the endpoint
of $\alpha _{j}$, for any $i,j$ we have 
\begin{equation*}
d_{\widehat{X}}(y_{i},y_{j})\leq 2(d_{X}(x_{j},x)+d_{X}(x_{i},x))
\end{equation*}%
and it follows from the convergence of $\{x_{j}\}$ that $\{y_{j}\}$ is
Cauchy. Moreover, since the inclusion of $\widehat{X}$ into $\widetilde{X}$
is uniformly continuous, $\{y_{i}\}$ is also Cauchy in $\widetilde{X}$. Now $%
\{y_{j}\}$ is actually a sequence in $\phi ^{-1}(x)$, which is an orbit of
the action of $\delta (X)$. But this orbit is complete (with the metric
induced by $\widetilde{X}$) by Proposition 36 of \cite{P} and therefore $%
y_{j}\rightarrow y\in \widetilde{X}$. Without loss of generality we may
suppose that $d(y_{j},y_{j+1})<2^{-j}$ for all $j$. Let $c_{j}$ be a path
joining $y_{j}$ and $y_{j+1}$ of length less than $2^{-j}$ in $\widehat{X}$.
We may parameterize $c_{1}$ on $[0,\frac{1}{2}]$, $c_{2}$ on $[\frac{1}{2},%
\frac{3}{4}]$, and so on. The concatenation of these paths is a distance
decreasing path from $[0,1)$ into $\widehat{X}$ and hence a uniformly
continuous path into $\widetilde{X}$. Since $y_{j}\rightarrow y$ in $%
\widetilde{X}$, this path has a unique continuous extension to a path $%
c:[0,1]\rightarrow \widetilde{X}$ from $y_{1}$ to $y$. Moreover, since $%
\widehat{\phi }$ is length preserving, $\widehat{\phi }\circ c$ is a path of
length at most $1$ in $X$. It now follows that $y\in \widehat{X}$ and $%
y_{j}\rightarrow y$ in $\widehat{X}$. By definition, $d(z_{j},y_{j})%
\rightarrow 0$ and hence $z_{j}\rightarrow y$ in $\widehat{X}$ as well.

For the fifth part, suppose that $g\in \delta (X)$ stabilizes $\widehat{X}$.
Then $\ast $ and $g(\ast )$ are joined by an path $\alpha $ such that $\phi
(\alpha )$ is rectifiable. Since $\phi (g(\ast ))=\phi (\ast )=\ast $, $\phi
(\alpha )$ is a loop representing an element $[\alpha ]\in \mu (X)$ with $%
\lambda ([\alpha ])=g$. On the other hand, suppose that $g\in \mu (X)$ and
let $\alpha $ be a rectifiable loop representing $g$. Then $\ast $ and $%
g(\ast )$ are joined by $\alpha _{L}$. Now if $x\in \widehat{X}$, we may
join $\ast $ to $x$ by a path $\beta $ such that $\phi \circ \beta $ is
rectifiable. Now $g\circ \beta $ joins $g(\ast )$ and $g(x)$, and $\phi
\circ g\circ \beta =g\circ \beta $ is rectifiable. Therefore the
concatenation $\gamma $ of $\alpha _{L}$ and $g\circ \beta $ is a path
joining $\ast $ and $g(x)$ such that $\phi \circ \gamma $ is rectifiable.
Therefore $g(x)\in \widehat{X}$ and we have shown that $\mu (X)$ is the
stabilizer of $\widehat{X}$. The fact that $\mu (X)$ acts isometrically
follows from the fact that for any path $\alpha $ in $\widehat{X}$ and $g\in
\mu (X)$, $\phi \circ \alpha =\phi \circ g\circ \alpha $. Since $\delta (X)$
acts freely, so does $\mu (X)$. Moreover, the orbits of $\delta (X)$ are
closed. Suppose $p\in \widehat{X}$. Since $\mu (X)$ is the stabilizer of $%
\widehat{X}$ in $\delta (X)$, the orbit $\mu (X)p$ is the intersection of
the orbit $\delta (X)p$ with $\widehat{X}$ and hence is closed in the
subspace topology of $\widehat{X}\subset \widetilde{X}$. But we already know
that the inclusion of $\widehat{X}$ with the length is continuous, hence the
orbits of $\mu (X)$ are closed in this topology. Since $\widehat{\phi }$ is
already a metric quotient, the final part of the fifth statement is simply
the observation that the orbits of $\mu (X)$ are exactly the pre-images of
points with respect to $\widehat{\phi }$.
\end{proof}

\begin{corollary}
\label{same}There is a unique (up to basepoint choice) URL-map $\theta :%
\overline{X}\rightarrow \widehat{X}$ such that $\overline{\phi }=\widehat{%
\phi }\circ \theta $. Moreover,

\begin{enumerate}
\item The restriction $\theta ^{\ast }:\Lambda (X)\rightarrow \mu (X)$ is a
surjective homomorphism.

\item If $X$ is uniformly $1$-dimensional then $\theta $ is an isometry and $%
\theta ^{\ast }$ is an isomorphism.
\end{enumerate}
\end{corollary}

\begin{proof}
The first part is a corollary of Theorems \ref{m2} and \ref{core}. The fact
that $\theta ^{\ast }$ is a homomorphism follows from the fact that the
cancelled concatenation\ of two paths is homotopic to the concatenation. Now
suppose that $\gamma \in \mu (X)$. Take a $\rho $-path $c$ joining $\ast $
and $\gamma (\ast )$. Then $\widehat{\phi }\circ c\in \Lambda (X)$ with $%
\theta ^{\ast }(\widehat{\phi }\circ \gamma )=\gamma $ by definition,
proving that $\theta ^{\ast }$ is surjective.

To prove the last part we need only show that $\theta $ is an injection when 
$X$ is uniformly $1$-dimensional. But if $\theta (c)=\theta (d)$ then the
lifts $c_{L}$ and $d_{L}$ must have the same endpoint. The compositions $%
c^{\prime },d^{\prime }$ of these paths with the inclusion of $\widehat{X}$
into the simply connected uniformly one-dimensional space $\widetilde{X}$
also have the same endpoints and hence must be fixed-endpoint homotopic by
Corollary \ref{runiq}. But then $c=\phi \circ c^{\prime }$ and $d=\phi \circ
d^{\prime }$ are also fixed-endpoint homotopic. By Corollary \ref{rep}, $c$
and $d$ must be identical.
\end{proof}

\begin{remark}
The kernel $K$ of $\theta ^{\ast }$ is of obvious interest, since $\pi
_{1}^{K}(X)=\Lambda (X)/K$ is isomorphic to $\mu (X)$. Now $K$ consists of
all elements of $\Lambda (X)$ that lift as loops via $\widehat{\phi }$, and
so $K=H_{L}$ if and only if $\widehat{\phi }$ is $H_{L}$-simply connected.
\end{remark}


\begin{example}
The infinite torus $\mathbb{T}^{\infty }=S^{1}\times S^{1}\times \cdot \cdot
\cdot $ can be metrized in the following natural way (cf. \cite{BPS}). Give
each $S^{1}$ a geodesic metric such that the diameters of the factors are
square summable, and apply the geodesic product metric to $\mathbb{T}%
^{\infty }$ (this is not the \textquotedblleft
topologist's\textquotedblright\ product metric, but the one that generalizes
the usual Euclidean or Riemannian product metrics, and requires square
summability of the diameters). For this particular metric we already
considered the metric core in \cite{BPS}, although we did not have a general
construction at the time and did not refer to it as the metric core. In fact
the metric core is simply separable Hilbert space $l^{2}$ consisting of all
square summable sequences inside $\widetilde{T^{\infty }}=\mathbb{R}\times 
\mathbb{R\times \cdot \cdot \cdot }$. Note that as is seen from the results
in \cite{BPS} there are non-rectifiable $1$-parameter subgroups in $\mathbb{T%
}^{\infty }$ that cannot be lifted to $l^{2}$ (although they do, of course,
lift to $\widetilde{T^{\infty }}$).
\end{example}

\section{Appendix}

The uniform universal covering (UU-covering) was defined in \cite{BPU}, and
provides an analog of the traditional universal covering for uniform spaces
that are not necessarily semi-locally simply connected or even locally
pathwise connected. In this appendix we will sketch the construction and
results of \cite{BPU} in the simplified setting of metric spaces, while also
establishing simplified notation and providing a reference for those
primarily interested in geodesic spaces as opposed to uniform spaces in full
generality.

Let $X$ be an arbitrary metric space. An $\varepsilon $-chain in $X$ is a
finite ordered sequence of points $\{x_{0},...,x_{n}\}$ such that for each $%
i $, $d(x_{i},x_{i+1})<\varepsilon $; the $\varepsilon $-chain is an $%
\varepsilon $-loop if $x_{0}=x_{n}$. Two $\varepsilon $-chains having the
same pair of endpoints are called $\varepsilon $-homotopic if one can be
obtained from the other via a finite number of steps, each of which involves
either adding or taking away a single point, maintaining an $\varepsilon $%
-chain at each step. The $\varepsilon $-homotopy equivalence class of an $%
\varepsilon $-chain $\gamma $ is denoted by $[\gamma ]_{\varepsilon }$. We
may choose any basepoint $\ast \in X$; all theorems that follow are
independent of the choice of basepoint, and when dealing with maps between
spaces we may always assume that the map takes the basepoint of one space to
the basepoint of the other.

The set of all $\varepsilon $-homotopy equivalence classes $[\gamma
]_{\varepsilon }$ of $\varepsilon $-chains starting at the basepoint $\ast $
is denoted by $X_{\varepsilon }$. We denote by $\delta _{\varepsilon }(X)$
the $\varepsilon $-deck group of $X$, which consists of $\varepsilon $%
-homotopy classes of $\varepsilon $-loops based at $\ast $, with the group
operation induced by concatenation, which also acts as a group of bijections
on $X_{\varepsilon }$ via concatenation. Moreover, $X$ is naturally
identified with the orbit space $\delta _{\varepsilon }(X)\backslash
X_{\varepsilon }$.

If $0<\delta <\varepsilon $ then every $\delta $-chain (resp. homotopy) may
be considered as an $\varepsilon $-chain (resp. homotopy) and therefore
there is a well-defined function $\phi _{\varepsilon \delta }:X_{\delta
}\rightarrow X_{\varepsilon }$ defined by $\phi _{\varepsilon \delta }\left( %
\left[ \gamma \right] _{\delta }\right) =\left[ \gamma \right] _{\varepsilon
}$. Since $\phi _{\varepsilon \delta }=\phi _{\varepsilon \alpha }\circ \phi
_{\alpha \delta }$ whenever $0<\delta <\alpha <\varepsilon $, we have an
inverse system $(X_{\varepsilon },\phi _{\varepsilon \delta })$ indexed by
the positive reals with the reverse ordering. The inverse limit of this
system is denoted by $\widetilde{X}$. Roughly speaking, elements of $%
\widetilde{X}$ may be thought of as collections of discrete homotopy
equivalence classes of finer and finer chains (rather than homotopy classes
of paths used to construct the traditional universal covering). The natural
projection $\phi :\widetilde{X}\rightarrow X$ takes an element $([\gamma
_{\varepsilon }]_{\varepsilon })$ of $\widetilde{X}$ to the common endpoint $%
x$ of all of the chains $\gamma _{\varepsilon }$. The set of $([\gamma
_{\varepsilon }]_{\varepsilon })\in \widetilde{X}$ such that each $\gamma
_{\varepsilon }$ is an $\varepsilon $-loop forms a group $\delta (X)$ with
respect to concatenation, called the \textit{uniform fundamental group. }(In 
\cite{BPU} we denoted this group by $\delta _{1}(X)$ and called it the
\textquotedblleft deck group\textquotedblright .)\textit{\ }$\delta (X)$
acts freely on $\widetilde{X}$ by concatenation.

Now suppose that $X$ is a length space. The length metric of $X$ may be
lifted to a length metic on $X_{\varepsilon }$ in such a way that the
bonding maps $\phi _{\varepsilon \delta }$ are local isometries and
traditional covering maps. The projection $\phi :\widetilde{X}\rightarrow X$
is surjective and is called the UU-covering of $X$, although it is not a
traditional covering map.

\begin{notation}
\label{indexnot}In any countable inverse limit construction, one may use a
cofinal sequence of indices to obtain the inverse limit. For simplicity we
will often use the sequence $2^{-i}$ to index the system for $i=1,2,...$,
and let $X:=X_{0}$, $X_{i}:=X_{2^{-i}}$, $\phi _{ij}:=\phi _{2^{-i}2^{-j}}$, 
$\phi _{i}:=\phi $, and $\delta _{i}(X):=\delta _{2^{-i}}(X)$. We will
denote elements of $\widetilde{X}$ as sequences $(x_{i})$ with $x_{i}\in
X_{i}$.
\end{notation}

While $\widetilde{X}$ is metrizable, there is no hope of finding a
compatible length on $\widetilde{X}$ because this space may not be pathwise
connected. However, there is a metric $\widetilde{d}$ on $\widetilde{X}$
compatible with the inverse limit uniform structure such that $\delta (X)$
acts isometrically and the orbits of $\delta (X)$ in $\widetilde{X}$ are
closed. Moreover, the orbit space $\delta (X)\backslash\widetilde{X}$ with
the quotent metric is uniformly homeomorphic to $X$ via the mapping $%
x\leftrightarrow \phi ^{-1}(x)$.

If $X$ is a compact geodesic space that is semilocally simply connected then
the fundamental inverse system stabilizes at some sufficiently large $i$,
and the projection $\phi _{\infty i}:X_{i}\rightarrow X$ is the traditional
universal covering of $X$, and $\delta (X)=\delta _{i}(X)$ is the
fundamental group of $X$.

The UU-covering $\phi :\widetilde{X}\rightarrow X$ has the following lifting
property. In \cite{BPU} we defined a notion of \textquotedblleft
universal\textquotedblright\ uniform space that is an analog of simply
connected in the traditional setting. If $Y$ is universal and $%
f:Y\rightarrow X$ is uniformly continuous, then there is a unique uniformly
continuous function $f_{L}:Y\rightarrow \widetilde{X}$ such that $\phi \circ
f_{L}=f$ and $f_{L}(\ast )=\ast $. The function $f_{L}$ will be called the 
\textit{lift} of $f$. Since real segments and their cartesian products are
universal, one may lift paths and homotopies to $\widetilde{X}$. There is a
natural homomorphism $\lambda :\pi _{1}(X)\rightarrow \delta (X)$ defined as
follows: If $\alpha $ is a loop representing an element of $\pi _{1}(X)$
based at $\ast $ then $\lambda ([\alpha ])$ is the unique element of $\delta
(X)$ that carries the basepoint $\ast $ to the endpoint of $\alpha _{L}$.
This homomorphism is injective if and only if the pathwise connected
component $P$ of $\ast $ in $\widetilde{X}$ is simply connected, and
surjective if and only if $\widetilde{X}$ is pathwise connected.

Although $\widetilde{X}$ need not be pathwise connected, from the results of 
\cite{BPU} we know that the pathwise connected component of $\widetilde{X}$
is dense in $\widetilde{X}$ and in particular $\widetilde{X}$ is connected.
Note that $\phi _{i}:\widetilde{X}\rightarrow X_{i}$ is the UU-covering of $%
X_{i}$ for any $i$.

The UU-covering also has an induced mapping property: If $f:X\rightarrow Y$
is uniformly continuous then there is a unique (basepoint preserving)
uniformly continuous function $\tilde{f}:\tilde{X}\rightarrow \tilde{Y}$
such that if $\phi :\widetilde{X}\rightarrow X$ and $\psi :\widetilde{Y}%
\rightarrow Y$ are the UU-coverings, $f\circ \phi =\psi \circ \widetilde{f}$%
. The lift is functorial and the restriction of $\widetilde{f}$ to $\delta
(X)\subset \widetilde{X}$ is a homomorphism $f_{\ast }:\delta (X)\rightarrow
\delta (Y)\subset \widetilde{Y}$.

Note that the mappings $\phi _{ij}:X_{j}\rightarrow X_{i}$ are by
construction URL-maps.

\begin{remark}
\label{light} The UU-covering map of any length space is also light; in fact
point pre-images are inverse limits of discrete spaces and so are totally
disconnected. Note that this statement, in particular, is also true in the
more general setting of coverable uniform spaces described in \cite{BPU}.
\end{remark}

We proved in Proposition 92 and Theorem 93 from \cite{BPU} the following:

\begin{proposition}
\label{inj} If $X$ is a coverable uniform uniformly one-dimensional space,
then $\tilde{X}$ is simply connected and contains no topological circles,
and the homomorphism $\lambda :\pi _{1}(X)\rightarrow \delta (X)$ is
injective.
\end{proposition}

\begin{remark}
As we showed in \cite{BPU}, every connected, locally arcwise connected
compact uniform space, hence any Peano continuum $X$, is coverable. But one
may also apply the Bing-Moise theorem to obtain a compatible geodesic metric
and use the construction described in this appendix.
\end{remark}

\end{document}